\pdfoutput=1
\RequirePackage{ifpdf}
\ifpdf 
\documentclass[pdftex]{sigma}
\else
\documentclass{sigma}
\fi

\usepackage{cases}
\usepackage{mathtools}

\numberwithin{equation}{section}

\newtheorem{Theorem}{Theorem}[section]
\newtheorem*{Theorem*}{Theorem}
\newtheorem*{TheoremMain}{Theorem~\ref{thm:main}}

\newtheorem{Lemma}[Theorem]{Lemma}
\newtheorem{Proposition}[Theorem]{Proposition}
\newtheorem{Claim}[Theorem]{Claim}

\theoremstyle{definition}
\newtheorem{Definition}[Theorem]{Definition}

\newtheorem*{Notation}{Notation}

\newtheorem{Remark}[Theorem]{Remark}

 \newcommand{\BbR}{{\mathbb R}}
 \newcommand{\BbP}{{\mathbb P}}
 \newcommand{\BbC}{{\mathbb C}}
 \newcommand{\BbZ}{{\mathbb Z}}
 \newcommand{\pdo}{\Psi{\rm DO}}
 
 \newcommand{\dvol}{{\rm dvol}}

 \newcommand{\diff}{{\rm Diff} }
 \newcommand{\maps}{{\rm Maps} }
 \newcommand{\kk}{2k-1 }

\newcommand{\Z}{\BbZ}
\newcommand{\R}{\BbR}
\newcommand{\C}{\BbC}
\newcommand{\Tr}{{\rm Tr}}
\newcommand{\tr}{{\rm tr}}
\newcommand{\calL}{\mathcal{L}}

 \newcommand{\resw}{{\rm res}^{\rm W}}

 \newcommand{\be}{\eta}
 \newcommand{\bo}{\overline\Omega}
 \newcommand{\bom}{\bo_{\overline M}}
 \newcommand{\bn}{\nabla}
 
 \newcommand{\bbn}{\overline{\nabla}}
 \newcommand{\bxi}{\xi}
 
\newcommand{\br}{ R}
\newcommand{\brr}{\overline R}
\newcommand{\bgg}{\bar g}

\newcommand{\bg}{ g}
\newcommand{\bx}{X}
\newcommand{\by}{Y}

\newcommand{\bm}{ M}
\newcommand{\ocs}{\widetilde{{\rm CS}}}
\newcommand{\tp}{\tilde p}
\newcommand{\xl}{X^L}
\newcommand{\yl}{Y^L}
\newcommand{\zl}{Z^L}

\newcommand{\sgn}{{\rm sgn}}
\newcommand{\dg}{\dot\gamma}
\newcommand{\ints}{\int_{S^1}}

\newcommand{\eb}{\overline e}

\newcommand{\om}{\overline{M}}
\newcommand{\omz}{M_0}
\newcommand{\omp}{M_p}

\newcommand{\tg}{\tilde g}
\newcommand{\barcurvature}[2]{R_{#1}{}^{#2}}

\newcommand{\ssuzero}{\sum_{\sigma_0 = 0} \sgn(\sigma)}

\newcommand{\isom}{{\rm Isom}}
\newcommand{\isombm}{\isom(\bm)}
\newcommand{\itbm}{[0,1]\times M}
\newcommand{\itm}{[0,1]\times M}
\newcommand{\ism}{\isom}

\begin{document}

\allowdisplaybreaks

\newcommand{\arXivNumber}{2011.01800}

\renewcommand{\thefootnote}{}

\renewcommand{\PaperNumber}{054}

\FirstPageHeading

\ShortArticleName{The Geometry of Loop Spaces III: Isometry Groups of Contact Manifolds}

\ArticleName{The Geometry of Loop Spaces III:\\ Isometry Groups of Contact Manifolds\footnote{This paper is a~contribution to the Special Issue on Interactions of Poisson Geometry, Lie Theory and Symmetry in honor of Rui Loja Fernandes for his 60th birthday. The~full collection is available at \href{https://sigma-journal.com/Fernandes.html}{https://sigma-journal.com/Fernandes.html}}}

\Author{Satoshi EGI~$^a$, Yoshiaki MAEDA~$^{\rm b}$ and Steven ROSENBERG~$^{\rm c}$}

\AuthorNameForHeading{S.~Egi, Y.~Maeda and S.~Rosenberg}

\Address{$^{\rm a)}$~Rakuten Institute of Technology, Rakuten, Inc., Japan}
\EmailD{\mail{satoshi.egi@rakuten.com}}
\URLaddressD{\url{https://www.egison.org/~egi/}}

\Address{$^{\rm b)}$~Tohoku Forum for Creativity, Tohoku University, Japan}
\EmailD{\mail{yoshimaeda@tohoku.ac.jp}}

\Address{$^{\rm c)}$~Boston University, USA}
\EmailD{\mail{sr@math.bu.edu}}
\URLaddressD{\url{https://math.bu.edu/people/sr}}

\ArticleDates{Received April 25, 2025, in final form May 03, 2026; Published online May 29, 2026}

\Abstract{Let $M_p$ be a circle bundle with first Chern class $p[\omega]$ over a closed $4n$-dimensional integral symplectic manifold \smash{$\bigl(\overline{M},\omega\bigr)$}. Equivalently, $M_p$ is a~closed contact $(4n+1)$-manifold whose Reeb orbits are all closed and have the same period. For a metric~$g$ on~$M_p$ compatible with the symplectic structure and the geometry of the circle fiber, we use Wodzicki--Chern--Simons forms on the loop space $LM_p$ to prove that $\pi_1({\rm Isom}(M_p,g))$ is infinite for ${|p| \gg 0}$. We~also give the first high-dimensional examples of nonvanishing Wodzicki--Pontryagin~forms.}


\Classification{35S99; 53D05; 53D10; 58D15}

\renewcommand{\thefootnote}{\arabic{footnote}}
\setcounter{footnote}{0}

\section{Introduction}

In this paper, we study
circle bundles $\omp$ with first Chern class $p[\omega]$ over a closed $4n$-dimensional integral
symplectic manifold \smash{$\bigl(\overline M,\omega\bigr)$}. Equivalently, $\omp$ is a
closed
$(4n+1)$-dimen\-sion\-al contact manifold with closed Reeb orbits \cite[Theorem~3.9]{Blair}. We give $\omp$ a Riemannian metric $\bg$ compatible with the symplectic structure on the base and the natural connection on the fibers. There is an infinite-dimensional family of such metrics.
We use Wodzicki--Chern--Simons (WCS) forms~\smash{$\ocs{}^{\rm W}_{2k-1}$} on the loop space $L\omp$ to determine that
 $\pi_1({\rm Isom}(\omp,\bg))$ is infinite for $|p| \gg 0$. This extends results for circle bundles over K\"ahler surfaces in \cite{MRT4} (as corrected in \cite{MRT4-cor}) to symplectic manifolds in arbitrarily high dimensions.

In general, the isometry group $\isom(\bm) = \isom(\bm,g)$ of a closed manifold $\bm$ is well-known to be a compact Lie group. It follows that $\isom(\bm)$ is isomorphic to
$\bigl(T^k\ltimes G\bigr)/F$, where $T^k$ is the $k$-torus, $G$ is a semisimple Lie group, and $F$ is a finite group \cite[Theorem~6.9]{Bredon}.
Since~$\pi_1(G)$ is finite, $ \pi_1(\isom(\bm)) :=\pi_1(\isom(\bm), {\rm Id}) $
 is infinite iff $k \geq 1$.
However, it seems difficult in general to determine $k$.

Suppose $M$ admits a nontrivial circle action $a\colon S^1\times
M\to M$ via isometries. This gives a~loop $a^I\colon S^1\to \isombm$ of isometries and hence an element of
 $\pi_1(\isombm)$. If the circle action is free, $M$ is the total space of a circle bundle over the orbit space $\om$, with the action given by rotation of the circle fibers. It is natural to conjecture that
 the class \smash{\raisebox{0.2pt}{$\bigl[a^I\bigr]\in \pi_1(\isombm)$}} has~infinite order. This is not always true: for the canonical bundle \raisebox{-0.2pt}{$ M= S^{2n+1}$} over $\overline M = \BbC\BbP^n$,
 the fiber rotation is an isometry of the standard metric on $S^{2n+1}$.
In fact, \smash{$\bigl[a^I\bigr]$} is the~generator of \smash{$\pi_1\bigl(\isom\bigl(S^{2n+1}\bigr)\bigr) = \pi_1({\rm SO}(2n+2))\simeq \Z_2$}.

In this example, the first Chern number of the canonical bundle is $1$.
The main result is that for sufficiently high Chern number, rotation in the circle fiber gives an element of infinite order in $\pi_1(\isombm)$. More precisely, we have the following.

\begin{TheoremMain}
Let $\bigl(\om,\omega\bigr)$ be a closed integral symplectic manifold of dimension $4n$. For $p\in \Z$, let $\omp$ be the circle bundle over $\om$ with first Chern class $p[\omega]$. Choose a metric $g$ on $\omp$ compatible with an almost complex structure on $\om$ as in \eqref{eq:gp}. Then for $|p| \gg 0 $,
 $\pi_1(\isom(\omp,g))$ is
infinite. Equivalently, let $M$ be a closed $(4n+1)$-dimensional contact manifold whose Reeb orbits are all closed and have the same period. Then $M$ covers infinitely many such contact manifolds $\omp$ with
$\pi_1({\rm Isom}(\omp,g))$~infinite.
\end{TheoremMain}

In the concrete example of $\BbC\BbP^{2}$, we proved in~\cite{MRT4} that $\pi_1(\isom(\omp))$ is infinite for
$p\neq \pm1$. In fact, the only example we know where~$a^I$ does not have infinite order in $\pi_1({\rm Isom}(\omp))$ is for~$\BbC\BbP^{n}$.

In Section~\ref{section2}, we give background material on pseudodifferential operators and WCS forms on
loop spaces. In Section~\ref{section3}, we prove the main result, both by direct calculation and computer verification.
The Stokes' theorem arguments (Propositions~\ref{customprop:app} and~\ref{prop:temp}) used in the proof depend on the key identity
\begin{equation} \label{eq:intro}
F^{L,*}{\rm d}_{LM}\ocs^{\rm W}_{2k-1} = 0\in \Lambda^{2k}([0,1]\times M),
\end{equation}
where $F^L\colon [0,1] \times M \to LM$ is induced by $F\colon [0,1]\times S^1\times M\to M$, a homotopy of loops of isometries. This identity probably fails for a homotopy of loops of diffeomorphisms, which led to errors in \cite{MRT4}, now corrected in Appendix~\ref{app:oldB}.

In Section~\ref{section4}, we apply our theory to the first example of a symplectic, non-K\"ahler manifold, due to Kodaira and Thurston. Through explicit calculations, we get the results in Theorem~\ref{thm:main} for all~$p$.

In Section~\ref{section5},
we introduce a second theme in the paper, the investigation of when WCS forms are closed. This would give interesting de Rham cohomology classes on $LM$. In finite dimensions, by (\ref{dcs}) a classical Chern--Simons form ${\rm CS}(\nabla_1, \nabla_0)$ on a bundle is closed if the connections $\nabla_0$, $\nabla_1$ are flat, or for dimension reasons. We had hoped to similarly find closed WCS forms on the tangent bundle to $LM$, where flatness or dimension restrictions are not available.
 On $LM$, \smash{$ {\rm d}_{LM}\ocs{}^{\rm W}_{2k-1}$} is a Wodzicki--Pontryagin (WP) form by~(\ref{dwcs}). Thus all WP
 classes are trivial, which is why we focus on secondary classes.

In Theorem~\ref{nonv}, we construct the first published examples of nonvanishing WP forms. The reasoning is as follows: if
\smash{$0= {\rm d}_{LM}\ocs{}^{\rm W}_{2k-1} \in \Lambda^{2k}(LM)$}, then (\ref{eq:intro}) holds for general families~$F$. (This vanishing was incorrectly assumed in~\cite{MRT4}.) In this case, the Stokes' Theorem argument and calculations in \cite[Section~3]{MRT4} apply in particular when $F$ is a homotopy through loops of smooth simple homotopy equivalences. Arguing as in Theorem~\ref{thm:main}, we would obtain that~$\pi_1(C^\infty{\rm Aut}(\mathcal{L}_p))$ is infinite, where $C^\infty{\rm Aut}(\mathcal{L}_p)$ is the space of smooth simple homotopy equivalences of a specific
five-dimensional lens space $\mathcal{L}_p$. This contradicts a result in~\cite{HJ}, so these WP forms are nontrivial. We note that these examples are $6$-forms; in finite dimensions, there are no degree~$6$ Pontryagin forms. The WP forms were predicted to vanish in \cite{lrst}, and \cite{alh-unpublished} gave the first counterexample.

Also in Section~\ref{section5}, we relate Pontryagin forms on $\om$ to WCS forms on
$M_p$, if $\om$ is K\"ahler (Proposition~\ref{app1}). If $\om$ is only symplectic, the proof breaks down, and no such relation is known.

Appendices~\ref{app:oldB},~\ref{app:B},~\ref{app:claim} and~\ref{app:E} give proofs of technical results. In Appendix~\ref{appb},
we discuss why symplectic manifolds of dimension $4n+2$ are more difficult to treat. The online files~\cite{Egi1, Egi3} include a particularly long calculation for the Kodaira--Thurston example and computer codes verifying the main results.

\begin{Notation}\quad
\begin{enumerate}\itemsep=0pt
\item[(i)] General odd-dimensional forms have degree
$2k-1$. We will work on specific circle bundles~$M_p$ over a symplectic manifold of dimension $4n$, so dim$(M_p) = 4n+1$. The WCS forms of interest are forms of degree $4n+1$ on $LM_p$.

\item[(ii)] Our conventions for the curvature tensor are as follows: for a Riemannian manifold $(M,g)$ with Levi-Civita connection $\nabla$, the curvature operator $R\colon TM^{\otimes 3} \to TM$ is defined by
\[R(A, B) C = \nabla_{A} \nabla_{B} C - \nabla_{B} \nabla_{A} C - \nabla_{[A,B]} C.\]
The components of the curvature tensor are given by
\begin{equation}\label{curvconv}\br(\partial_k, \partial_j)\partial_b = \br_{kjb}^{\ \ \ a}\partial_a,\qquad
\br_{kjba} = g(\br(\partial_k,\partial_j)\partial_b,\partial_a) =
\langle \br(\partial_k,\partial_j)\partial_b,\partial_a\rangle.
\end{equation}
In local coordinates, the matrix of curvature two-forms $\Omega \!=\! \Omega_g\in \Lambda^2(M, {\rm Hom}(TM,TM))$~is%
\begin{equation}\label{eq:cc}\Omega(\partial_k,\partial_j)_b{}^a = R_{jkb}{}^a.
\end{equation}
In these formulas and throughout the paper, we use Einstein summation convention.
\end{enumerate}
\end{Notation}

\section{Background material}\label{section2}
\subsection{Finite-dimensional background material}
The complexified tangent bundle
of a Riemannian manifold $\bigl(M^{4n},g\bigr)$ has Chern character
${\rm ch}(M) \in H^{{\rm ev}}(M,\R)$
with $2k$-component
\begin{equation}\label{onea}{\rm ch}_{[2k]}(M) = \frac{1}{k!(2\pi)^k}\bigl[\Tr\bigl(\Omega^k\bigr)\bigr]\in H^{2k}(M,\R),
\end{equation}
where $\Omega = \Omega_M$ is the curvature form of $g$. There are associated Pontryagin-type forms
$\tp_k(\Omega) = \smash{(-1)^k/\bigl[(2k)!(2\pi)^{2k}\bigr]\Tr\bigl(\Omega^{2k}\bigr)}$ and classes
\[\tp_{k}(M)= [ \tp_k(\Omega)] = (-1)^k {\rm ch}_{[2k]}(M)\in H^{4k}(M,\R).\]

The usual Pontryagin classes $p_k(M) := (-1)^k c_{2k}(M)$ are built from the even Chern classes $c_{2k}(M)$.
By invariant theory for ${\rm SO}(n)$, the rings generated by $\{\tp_{k}\}$ and
$\{p_{k}\}$ are the same; this reduces to Newton's identities relating the elementary symmetric functions in $\lambda_1,
\dots ,\lambda_n$ to~$\sum \lambda_i,\dots,\sum \lambda_i^n$ \cite[Section~16]{milnor}.

As part of Chern--Weil theory, for connections $\nabla^0$, $\nabla^1$ on $TM$
with curvature forms $\Omega^0$, $\Omega^1$, the Chern--Simons form
\smash{${\rm CS}_{4k-1}\bigl(\nabla^0,\nabla^1\bigr)
\in \Lambda^{4k-1}(M)$},
\[{\rm CS}_{4k-1}\bigl(\nabla^0,\nabla^1\bigr) = 2k\int_0^1 \Tr((\omega_1 - \omega_0)\wedge \overbrace{\Omega_t
\wedge\cdots\wedge\Omega_t}^{2k-1}) \, {\rm d}t,\]
 satisfies
\begin{equation}\label{dcs} {\rm dCS}_{4k-1}\bigl(\nabla^0,\nabla^1\bigr) = \tp_k\bigl(\Omega^0\bigr) - \tp_k\bigl(\Omega^1\bigr).
\end{equation}
Here $\omega_t =t\omega_0+(1-t)\omega_1$, $\Omega_t = {\rm d}\omega_t+\omega_t\wedge\omega_t$.

\subsection{Infinite-dimensional background material}
This material is taken from \cite{MRT3, MRT4}. Let $(M,g)$ be a Riemannian manifold. For fixed $s\gg 0$, the loop space $LM$ of $s$-differentiable loops is a Banach manifold with tangent space at a loop $\gamma\colon S^1\to M$ given by $T_\gamma LM =\Gamma\bigl(\gamma^*TM\to S^1\bigr)$, where the sections of the pullback bundle are
$s$-differentiable. $LM$
has two preferred connections, the $L^2$ or $s=0$ Levi-Civita connection $\nabla^0$ associated to the $L^2$ inner product $\langle\, ,\, \rangle_0$, and the $s=1$ Levi-Civita connection $\nabla^1$ associated to the inner products
$\langle\, ,\, \rangle_1$:
\[\langle X,Y\rangle_0 = \int_{S^1} g(X_t,Y_t)_{\gamma(\theta)}\, {\rm d}t,\qquad
\langle X,Y\rangle_1 = \int_{S^1} g((1+ \Delta)X_t,Y_t)_{\gamma(\theta)}\, {\rm d}t.\]
Here
$\Delta = \nabla^*\nabla$ is the Laplacian associated to the pullback connection $\nabla = \gamma^*\nabla^M$ of the Levi-Civita connection $\nabla^M$ on $M$. While the connection and curvature forms for $\nabla^0$
at $\gamma$ take values in End$(\gamma^*TM)$, the corresponding forms for $\nabla^1$ take values in
$\pdo_{\leq 0}$, the Lie algebra of zeroth
order
pseudodifferential operators ($\pdo$s) on $\Gamma(\gamma^*TM\otimes \C)$, with the understanding that zeroth order means order at most zero. Since endomorphisms of a bundle are zeroth order $\pdo$s, we can consider $\nabla^0$, $\nabla^1$ to be $\pdo_0^*$-connections, where the Lie group
$\pdo_0^*$ of zeroth order invertible $\pdo$s with bounded inverse has Lie algebra $\pdo_{\leq 0}$.
 In particular, the curvature forms for these connections take values in $\pdo_{\leq 0}$.

In contrast to finite dimensions, there are two natural traces on $\pdo_{\leq 0}$. Recall that a zeroth order $\pdo$
$P$ on $\Gamma(\gamma^*TM\otimes \C)$ has a symbol sequence \smash{$P\sim \sum_{k=0}^\infty \sigma^P_{-k}(x,\xi)$}, where
$x\in S^1$, ${\xi\in T^*_xS^1}$; for $\pi\colon T^*S^1\to S^1$ the projection, $\sigma^P_{-k}(x,\xi)\in {\rm End}\bigl(\pi^*\gamma^*TM|_{(x,\xi)}\bigr)$ is homogeneous of degree $-k$ in $\xi$. The first trace is the leading order trace
\[\Tr^{{\rm lo}}(P) =\frac{1}{4\pi} \int_{S^*S^1} \tr (\sigma_0(x,\xi) )\, {\rm d}\xi{\rm d}x,\]
where $S^*S^1$ is the unit cotangent bundle of $S^1$. For example, if $P\in {\rm End}(\gamma^*TM\otimes \C)$, then
$\Tr^{\rm lo}(P) = (1/2)\int_{S^1} \tr(P(x)) {\rm d}x$. The second
 is the Wodzicki residue (see \cite{fgls})
 \[\resw(P) = \frac{1}{4\pi}\int_{S^*S^1} \tr ( \sigma_{-1}(x,\xi)) \, {\rm d}\xi{\rm d}x.\]
 For $P\in {\rm End}(\gamma^*TM\otimes \C)$, $\resw(P) = 0$. The trace in (\ref{onea}) can be replaced by either trace to give a theory of characteristic classes on $TLM$:
 \begin{gather*}{\rm ch}_{[2k]}^{\rm lo}(LM) := \frac{1}{k!}\bigl[\Tr^{\rm lo}\bigl(\Omega^k\bigr)\bigr]\in H^{2k}(LM,\R), \\
 {\rm ch}_{[2k]}^{\rm W}(LM) := \frac{1}{k!}\bigl[\resw\bigl(\Omega^k\bigr)\bigr]\in H^{2k}(LM,\R).
 \end{gather*}
 In fact, the Chern character \smash{${\rm ch}_{[2k]}^{\rm W}(LM)$} always vanish, while there are many examples of nonvanishing \smash{${\rm ch}_{[2k]}^{\rm lo}(LM)$} \cite{lrst}.
 In this paper, we only consider the Wodzicki residue trace.
There are corresponding Wodzicki--Pontryagin (WP) classes
\[p_k^{\rm W}(LM),\tp_k^{\rm W}(LM) \in H^{4k}(LM,\R).\]
Since these classes vanish, we focus on the associated Wodzicki--Chern--Simons (WCS) forms
\[\ocs^{\rm W}_{2k-1} = k
\int_0^1 \resw((\omega_1 - \omega_0)\wedge \overbrace{\Omega_t
\wedge\cdots\wedge\Omega_t}^{k-1}) \, {\rm d}t\in\Lambda^{2k-1}(LM).\]
Fix a loop $\gamma(\theta)\in LM$ and complexified tangent vectors $X_1, \dots,X_{2k-1}\in \Gamma(\gamma^*TM\otimes \C)$ at $\gamma$.
By \cite[equation~(2.9)]{MRT4}, for the $L^2$ and $s=1$ Sobolev connections, we have
 \begin{align}
&\ocs^{\rm W}_{2k-1,\gamma} (X_1,\dots ,X_{2k-1}) \nonumber\\
&\qquad{}=
\frac{k}{2^{k-2}}
\sum_{\sigma\in \mathfrak S_{2k-1}} {\rm sgn}(\sigma) \int_{S^1}\tr\bigl[
 \bigl(R\bigl(X_{\sigma(1)},\cdot\bigr)\dg\bigr)
 \Omega^{k-1}\bigl(X_{\sigma(2)},\dots,X_{\sigma(2k-1)} \bigr)\bigr]\nonumber\\
 &\qquad{}= \frac{k}{2^{k-2}}
\int_{S^1}
 K_{\nu \lambda_1 \dots \lambda_{2k-1}} (\gamma (\theta))
{\dot \gamma}^{\nu} (\theta )
X_{\sigma(1)}^{\lambda_{1}}(\theta ) \cdots X_{ \sigma(2k-1)}^{\lambda_{2k-1}}(\theta)
\, {\rm d}\theta,\label{csg}
\end{align}
\noindent where $R = R_g$, $\Omega = \Omega_g$ are the curvature tensor and curvature two-form of $g$,
 $\mathfrak S_{2k-1}$ is the permutation group of $\{1,\dots,2k-1\}$,
and by (\ref{eq:cc}),
\begin{align}\label{K_tensor}
& (2k-1)!K_{\nu \lambda_1 \dots \lambda_{2k-1}} \\
&\qquad{}= \sum_{\sigma\in \mathfrak S_{2k-1} } \sgn({\sigma})
R_{\lambda_{\sigma (1)} e_1 \nu}{}^{e_{2}}
R_{\lambda_{\sigma (2) } \lambda_{\sigma (3)} e_{3}}{}^{e_1}
R_{\lambda_{\sigma (4) } \lambda_{\sigma (5)} e_{4}}{}^{e_{3}}
\cdots
R_{\lambda_{\sigma (2k-2) } \lambda_{\sigma (2k-1)} e_{2}}{}^{e_{k-1}}.
\nonumber
\end{align}

The analog of (\ref{dcs}) in this context is
\smash{${\rm d}\ocs^{\rm W}_{4k-1}\bigl(\nabla^0,\nabla^1\bigr) = \tp^{\rm W}_k\bigl(\Omega^0\bigr) - \tp^{\rm W}_k\bigl(\Omega^1\bigr)$},
where $\Omega_0$, resp.~$\Omega_1$, are the curvature of the $L^2$, resp.~Sobolev $s=1$, metrics on $LM$.
Since $\Omega_0$ takes values in endomorphisms of $TLM$, its Wodzicki residue vanishes. Thus
\begin{equation}\label{dwcs}
{\rm d}\ocs^{\rm W}_{4k-1}\bigl(\nabla^0,\nabla^1\bigr) = - \tp_k^{\rm W}\bigl(\Omega^1\bigr)\in \Lambda^{4k}(LM).
\end{equation}
This proves that \smash{$\bigl[\tp_k^{\rm W}\bigl(\Omega^1\bigr)\bigr]$} is zero in $H^{4k}(LM)$, which of course is not true in general in finite dimensions. These forms are not necessarily pointwise zero (Theorem~\ref{nonv}), a new result.

\section{WCS forms for circle bundles over symplectic manifolds}\label{section3}

In Section~\ref{section3.1}, we prove the main result Theorem~\ref{thm:main}. We first discuss the Riemannian geometry of circle bundles $\omp$, $p\in \Z$, over symplectic manifolds $\bigl(\om,\omega\bigr)$, where $c_1(\omp) = (2\pi)^{-1}[p\omega]$. We compute the curvature as a function of $p$. Using the curvature calculations, we prove that the WCS class on $L\omp$ is a polynomial in $p^2$ with nonzero top coefficient. As we explain, this proves the Theorem. In Section~\ref{section3.2}, we discuss computer calculations that verify our calculations.

\subsection{Geometry of line bundles over integral symplectic manifolds}\label{section3.1}

Let $\bigl(\om,\omega\bigr)$ be a compact integral symplectic manifold of real dimension $4n$;
equivalently, $\om$ is projective algebraic. The symplectic form $\omega\in H_2\bigl(\om,\Z\bigr)$ determines a Riemannian metric $\bar g(X,Y) = -\omega(JX,Y)$, where $J$ is a compatible almost complex structure. The set of such $J$ is infinite-dimensional, so we obtain an infinite-dimensional family of metrics.

We compute the symplectic volume form of $\om$ in local coordinates, where
 $J = J_{ i}^{\ j} {\rm d}x^i \otimes \frac{\partial}{\partial x^j}$
 and $\omega = \omega_{ij} {\rm d}x^i\wedge {\rm d}x^j$.

 \begin{Lemma} \label{lem:Y.one}\quad
 \begin{enumerate}\itemsep=0pt
\item[$(i)$] $J_{ij} := \bgg_{bj}J^{\ b}_{ i} = \omega_{ij}$.

\item[$(ii)$]
The symplectic volume form of $\om$ is
\[\frac{1}{(2n)!}\omega^{2n} = \frac{1}{(2n)!2^{2n}}\sum_{\sigma\in \mathfrak S_{4n}}\sgn (\sigma) J_{\sigma(1)\sigma(2)} \cdots
 J_{\sigma(4n-1)\sigma(4n)}\, {\rm d}x^1\wedge\cdots\wedge {\rm d}x^{4n},\]
where $\mathfrak S_{4n}$ is the permutation group on $\{1,\dots,4n\}$.
\end{enumerate}
\end{Lemma}

Our convention throughout the paper is that $J_j^b = J_j^{\ b}$, so $\bgg_{ab}J_j^b = J_{ja}$, not $J_{aj}$.

\begin{proof} (i) Since $\omega(X,Y) = \bgg(JX,Y)$, we get
\[\omega_{ij} = \bgg(J\partial_i, \partial_j) = \bgg(J_i{}^b\partial_b, \partial_j) = J_i{}^b \bgg_{bj} = J_{ij}.\]

(ii) This follows from (i), since
\[\omega^{2n} =\frac{1}{2^{2n}}
\sum_{\sigma\in \mathfrak S_{4n}}\sgn (\sigma) \omega_{\sigma(1)\sigma(2)} \cdots
 \omega_{\sigma(4n-1)\sigma(4n)}\, {\rm d}x^1\wedge\cdots\wedge {\rm d}x^{4n}.\tag*{\qed}\]
\renewcommand{\qed}{}
\end{proof}

Because $\omega$ is integral, it has an associated line bundle $L = L_1$ over $\om$. Let $\omp$ be the total space of the circle bundle \smash{$L_p\stackrel{\pi}{\to} \om$} associated to $p\omega$ for $p\in\Z$.
$L_p$ comes with a connection $\eta= \eta_p$ with ${\rm d}\eta = p\pi^*\omega$, the curvature of $\eta$. The metric $\bgg$ induces a metric $g=g_p$ on $\omp$ by
\begin{equation}\label{eq:gp}\bg(\bx,\by)= \bg_p(X,Y) = \bar g(\pi_*\bx, \pi_*\by) + \be_p(\bx)_p\be(\by)
\end{equation}
(see \cite[p.~37]{Blair}).
We also denote $\bg(\bx,\by)$ by $\langle \bx, \by\rangle$.

Let $\xi$ be
a vector tangent to the circle fiber with $\eta( \xi) = 1$, and let $X^L$ denote the horizontal lift to $\omp$ of a tangent vector $ X$ to $\om$. We have $\eta\bigl(\xl\bigr) = 0$.

We compute the Levi-Civita connection $\bn$ for $\bg$ in terms of the Levi-Civita connection $\overline{\nabla}$ for~$\om$.

\begin{Lemma}\label{lem nabla}\quad
 \begin{enumerate}\itemsep=0pt
\item[$(i)$] $\bn_{\bxi}\bxi = {\mathcal L}_{\bxi} \bxi = 0$,

\item[$(ii)$]
\smash{$\bn_{\xl}\yl = \bigl(\overline{\nabla}_XY\bigr)^L -p\bgg(JX,Y)\bxi$,}

\item[$(iii)$] $ \bn_{\xl}\bxi = \bn_{\bxi}X^L
 = p(JX)^L$.
\end{enumerate}
 \end{Lemma}

 Here $\mathcal L$ is the Lie derivative.

\begin{proof}
(i) As in \cite[Section~3.2]{MRT4}, each circle fiber is the orbit of an isometric $S^1$ action on $\omp$, so each circle is a geodesic ($\bn_{\bxi}\bxi=0$), with $\bxi$ preserved by the action (${\mathcal L}_{\bxi} \bxi = 0$).
Alternatively, for the first part, since ${\rm d}\eta(\cdot, \bxi) = 0$, we get
$\mathcal L_{\bxi}\eta = {\rm d}i_{\bxi}\eta + i_{\bxi} {\rm d}\eta = {\rm d} 1 + {\rm d}\eta(\bxi,\cdot) = 0$.
Thus
$\mathcal L_{\bxi} \bg = \mathcal L_{\bxi} (\pi^*\bgg + \eta\otimes\eta) =0$, so $\bxi$ is a Killing vector field.
This implies $\bg(\bn_{\bxi}\bxi, Z) + \bg(\bxi, \bn_{\bxi}Z)=0$. Setting
 $Z = \bxi$ and then $Z\perp\bxi$, we get $\bn_{\bxi}\bxi = 0$.

(ii) We define $H(X,Y)\in \R, FX = F(X)\in T\om$ by
\begin{align}&\bn_{\xl} \yl= \bigl(\overline{\nabla}_XY\bigr)^L + H(X,Y)\bxi,\label{one a}\\
&\bn_{\xl} \bxi= (FX)^L,
\label{one b}
\end{align}
These definitions are valid, since for (\ref{one a}), it follows from \cite[Lemma 1]{oneill} that
$\pi_*\bigl(\bn_{\xl}\yl\bigr) = \overline{\nabla}_XY$, so
\smash{$\nabla_{X^{L}} Y^{L} = \bigl(\overline\nabla_{X} Y\bigr)^{L} + H(X, Y)\xi$} for some $H(X,Y)$.
For (\ref{one b}),
$\langle\xi, \xi\rangle=1$ implies
$\langle\nabla_{X^{L}}\xi, \xi\rangle=0$,
so $\bn_{\xl} \bxi = (FX)^L$ for some $FX$.

We note that $H(X,Y) = -H(Y,X)$: using $\eta(\bx) = \bg(\bxi,\bx)$, we get
\begin{align*} 0&= (\mathcal L_{\bxi}g) (\bx,\by) = (\bn_{\bx}\eta)(\by) +
(\bn_{\by}\eta)(\bx))\\
&= \bg(\bxi, \bn_{\bx}\by) + \bg(\bxi, \bn_{\by}\bx = H(X,Y) + H(Y,X).
\end{align*}
Thus
\begin{align*} p\omega(X,Y) &= {\rm d} \eta\bigl(X^{L}, Y^{L}\bigr)
= \frac{1}{2}\bigl(\nabla_{X^{L}} \eta\bigl(Y^{L}\bigr) - \eta\bigl(\bn_{X^{L}}Y^{L}\bigr) - \nabla_{Y^{L}}\bar \eta\bigl(X^{L}\bigr) - \eta\bigl(\nabla_{Y^{L}}X^{L}\bigr) \bigr) \nonumber\\
 &= \frac{1}{2}\bigl(- \eta\bigl(\nabla_{X^{L}}Y^L\bigr) + \eta\bigl(\nabla_{Y^{L}}X^{L}\bigr) \bigr) = -\bg\bigl(\bxi,\nabla_{X^{L}}Y^L\bigr) +
 \bg\bigl(\bxi,\nabla_{Y^{L}}X^L\bigr) \nonumber\\
&= \frac{1}{2}(-H(X,Y) + H(Y,X))
= -H(X,Y).
 \end{align*}
 This implies
 \[H(X,Y) = -p\omega(X,Y) = p\omega\bigl(J^2X,Y\bigr) = -pg(JX,Y).\]

(iii)
From $\mathcal{L}_{\xi} X^{L}=0$, we get the first equality in (iii):
\[\nabla_{\xi} X^{L} - \nabla_{X^{L}} \xi = \bigl[\xi, X^{L}\bigr] =
\mathcal{L}_{\xi} X^{L}=0.\]
(This also gives an alternative proof of (i): since $\bigl\langle \xl,\bxi\bigr\rangle = 0$, we have
\[\bigl\langle \bn_{\bxi}\bxi, \xl\bigr\rangle = \bigl\langle \bxi, \bn_{\bxi} \xl\bigr\rangle = \langle \bxi, \bn_{\xl} \bxi \rangle = 0.\]
Since $\langle\xi, \xi\rangle=1$ implies
\smash{$ \langle\nabla_{\xi}\xi, \xi\rangle=0$,} we get
\smash{$\bn_{\bxi}\bxi = 0$}. Another proof that the circle fibers are geodesics is in
\cite[Theorem 5.2.13]{Klingenberg}.)

It follows from $\bigl\langle Y^{L}, \xi\bigr\rangle=0$ that
\[
\bigl\langle\nabla_{X^{L}}Y^{L}, \xi\bigr\rangle + \bigl\langle Y^{L}, \nabla_{X^{L}}\xi\bigr\rangle = 0 \qquad \text{or}\qquad
H(X,Y) + \bigl\langle Y^{L}, (F X)^{L}\bigr\rangle=0.\]
Since $\bg\bigl(X^L, \yl\bigr) = \bgg(X,Y)$, we have
$-pg(JX,Y) = H(X,Y) = -g(FX, Y)$, so $FX = pJX$.
\end{proof}

The curvature tensor $\br$ of $\bg$ is related to the curvature tensor $\brr$ of $\bar g$ as follows.

\begin{Lemma} \label{lem1}
\begin{align*} {\rm (i)}\ & \bg\bigl(\br\bigl(X^L, Y^L\bigr)\zl,W^L\bigr) =
 \bgg\bigl( \brr(X,Y)Z,W\bigr) + p^2[-\bgg( JY,Z)\bgg( JX,W)\\
&\hphantom{ \bg\bigl(\br\bigl(X^L, Y^L\bigr)\zl,W^L\bigr) =}{}
+\bgg( JX,Z)\bgg( JY,W) +2\bgg( JX,Y)\bgg( JZ,W)],\\
{\rm (ii)}\ & \bg\bigl(\br\bigl(\xl, \yl\bigr)\zl, \bxi\bigr) = - p \bgg\bigl(\bigl(\bbn_{X} J\bigr) Y, Z\bigr)
+ p \bgg\bigl(\bigl(\bbn_{Y} J\bigr) X, Z\bigr),\\
{\rm (iii)}\ & \bg\bigl(\br\bigl(\xl,\bxi\bigr)\yl,\bxi\bigr) = -p^2\bgg(X,Y),\\
{\rm (iv)}\ & \bg\bigl(\br\bigl(\xl, \bxi\bigr)\yl, \zl\bigr) = p\bgg\bigl(\bigl(\bbn_XJ\bigr)Y,Z\bigr).
\end{align*}
\end{Lemma}

\begin{proof} (i) and (ii).
We have
\begin{gather*}\nabla_{X^{L}} \nabla_{Y^{L}} Z^{L}
 = \nabla_{X^{L}} \bigl( \bigl(\bbn_{Y} Z\bigr)^{L} - p \bgg (J Y, Z) \xi\bigr)\nonumber\\
 \hphantom{\nabla_{X^{L}} \nabla_{Y^{L}} Z^{L}}{}
 = \nabla_{X^{L}} \bigl(\bbn_{Y} Z\bigr)^{L} - p \nabla_{X^{L}} (\bgg (J Y, Z) \xi)\nonumber\\
 \hphantom{\nabla_{X^{L}} \nabla_{Y^{L}} Z^{L}}{}
 = \bbn_{X}\bigl(\bbn_{Y} Z\bigr)^{L} - p \bgg \bigl(J X, \bbn_{Y} Z \bigr) \xi - p \bigl[ X^{L} (\bgg(J Y, Z)) \xi + \bgg(J Y, Z) \nabla_{X^ L} \xi \bigr]\nonumber\\
 \hphantom{\nabla_{X^{L}} \nabla_{Y^{L}} Z^{L}}{}
 =\bigl(\bbn_{X} \bbn_{Y} Z\bigr)^{L} - p \bgg\bigl(J X, \bbn_{Y} Z\bigr) \xi
 - \bgg\bigl(\bigl(\bbn_{X} J\bigr) Y, Z\bigr) \xi + \bgg\bigl(J \bbn_{X} Y, Z\bigr) \xi \nonumber\\
 \hphantom{\nabla_{X^{L}} \nabla_{Y^{L}} Z^{L}=}{}
 + \bgg \bigl(J Y, \bbn_{X} Z\bigr) \xi + \bgg (JY, Z) \bigl(p (J X)^{L}\bigr),\nonumber\\
\bigl[X^{L}, Y^{L}\bigr]
 = \nabla_{X^ L} Y^{L} - \nabla_{Y^ L} X^{L} \nonumber\\
 \hphantom{\bigl[X^{L}, Y^{L}\bigr]}{}
 = \bigl(\bbn_{X} Y\bigr)^{L} - p \bgg(J X, Y) \xi - \bigl(\bbn_{Y} X\bigr)^{L} + p \bgg(J Y, X) \xi \nonumber\\
 \hphantom{\bigl[X^{L}, Y^{L}\bigr]}{}
 = -2 p \bgg(J X, Y) \xi + [X, Y]^{L}, \nonumber
\end{gather*}
so
\begin{align*}&R\bigl(X^{L}, Y^{L}\bigr) Z^{L}
 = \bigl(\bbn_{X} \bbn_{Y} Z\bigr)^{L} - p \bgg\bigl(J X, \bbn_{Y} Z\bigr) \xi \nonumber\\
 &\qquad\quad{} - p \bigl(\bgg\bigl(\bigl(\bbn_{X} J\bigr) Y, Z\bigr) \xi + \bgg\bigl(J \bbn_{X} Y, Z\bigr) \xi + \bgg\bigl(J Y, \bbn_{X} Z\bigr) \xi + p \bgg(J Y, Z) (J X)^{L}\bigr) \nonumber\\
 &\qquad\quad{} - \bigl[\bigl(\bbn_{Y} \bbn_{X} Z\bigr)^{L} - p \bgg\bigl(JY, \bbn_{X} Z\bigr) \xi \nonumber\\
 &\qquad\quad\hphantom{- \bigl[\, }{} - p \bigl(\bgg\bigl(\bigl(\bbn_{Y} J\bigr) X, Z\bigr) \xi + \bgg\bigl(J \bbn_{Y} X, Z\bigr) \xi + \bgg\bigl(J X, \bbn_{Y} Z\bigr) \xi + p \bgg(J X, Z) (J Y)^{L}\bigr)\bigr] \nonumber\\
 &\qquad\quad{} - \bigl[\bigl(\nabla_{[X,Y]^{L}} Z^{L}\bigr) - 2 p \bgg (J X, Y) \nabla_{\xi} Z^{L}\bigr] \nonumber\\
 &\qquad{}= (R(X, Y)Z)^L - p^2 \bgg(J Y, Z)(J X)^{L} + p^2 \bgg(J X, Z)(J Y)^{L} + 2 p^2 \bgg(JX, Y)(J Z)^{L} \nonumber\\
 &\qquad\quad{} - p \bgg\bigl(\bigl(\bbn_{X} J\bigr)Y, Z\bigr) \xi + p \bgg\bigl(\bigl(\bbn_{Y} J\bigr)X, Z\bigr) \xi. \nonumber
\end{align*}
Thus,
\begin{gather*}
g\bigl(R\bigl(X^{L}, Y^{L}\bigr)Z^{L}, W^{L}\bigr)
 = \bgg\bigl( \brr(X, Y) Z, W\bigr) - p^{2} \bgg(J Y, Z) \bgg(J X, W)\\
 \hphantom{g\bigl(R\bigl(X^{L}, Y^{L}\bigr)Z^{L}, W^{L}\bigr)=}{}
 + p^{2} \bgg(J X, Z) \bgg(J Y, W) + 2 p^{2} \bgg(J X, Y) \bgg(J Z, W),\\
\bg\bigl( R\bigl(X^{L}, Y^{L}\bigr)Z^{L}, \xi \bigr)
 = - p \bgg\bigl(\bigl(\bbn_{X} J\bigr) Y, Z\bigr) + p \bgg\bigl(\bigl(\bbn_{Y} J\bigr) X, Z\bigr).
\end{gather*}

(iii) and (iv). Using $\bigl[\xl,\bxi\bigr] = 0$ and Lemma~\ref{lem nabla}, we have
\begin{align*}
 \br\big(\xl,\bxi\big)\yl &= \bn_{\xl}\bn_{\bxi}\yl - \bn_{\bxi}\bn_{\xl}\yl\\
 & = \bn_{\xl}\bigl(p(JY)^L\bigr) - \bn_{\bxi}\bigl(\bigl(\overline{\nabla}_XY\bigr)^L - p\bgg(JX,Y)\bxi\bigr)\\
 &= p \bn_{\xl}\bigl((JY)^L\bigr)-p \bigl(J\bbn_XY\bigr)^L\\
 &= p\bigl[\bigl(\bbn_X(JY)\bigr)^L -p\bgg(JX,JY)\bxi\bigr]-p \bigl(J\bbn_XY\bigr)^L\\
 &= p\bigl(\bigl(\bbn_XJ\bigr)Y\bigr)^L - p^2 \bgg(X,Y)\bxi.
 \end{align*}
In other words,
\begin{align*}
 &\bg\bigl(\bigl( \br\bigl(\xl,\bxi\bigr)\yl),\zl\bigr)= p\bgg\bigl(\bigl(\bbn_XJ\bigr)Y, Z\bigr),\\
& \bg\bigl(\bigl( \br\bigl(\xl,\bxi\bigr)\yl\bigr),\bxi\bigr)= -p^2\bgg(X,Y).
\tag*{\qed}
\end{align*}\renewcommand{\qed}{}
 \end{proof}

 In fact, (ii) and (iv) are equivalent; this uses the symmetry of $\br$ and
 \begin{align*} 0 &= -{\rm d}\omega(X,Y,Z) = {\rm d}(\bgg(J\cdot,\cdot))(X,Y,Z)\\
 & \Rightarrow
 \bgg\bigl(\bigl(\bbn_X J\bigr) Y,Z\bigr) + \bgg\bigl(\bigl(\bbn_Y J\bigr)Z,X\bigr) - \bgg\bigl(\bigl(\bbn_Z J\bigr)Y,X\bigr) = 0.
 \end{align*}

 Here is the main result.

\begin{Theorem} \label{thm:main} Let $\bigl(\om,\omega\bigr)$ be a closed integral symplectic manifold of dimension $4n$.
Then for $|p| \gg 0 $,
 $\pi_1({\rm Isom}(\omp,g_p))$ is infinite. Equivalently, let $M$ be a closed $(4n+1)$-dimensional contact manifold whose Reeb orbits are all closed and have the same period. Then $M$ covers infinitely many such contact manifolds~$\omp$ with
$\pi_1({\rm Isom}(\omp,g_p))$ infinite.
\end{Theorem}

In fact, $\omp$ is diffeomorphic to $M_{-p}$, since $L_p$ is diffeomorphic to $L_{-p} = L_p^*$
via the fiberwise map $v\mapsto \langle \cdot, v\rangle$.

 For the equivalence, we note that the line bundle $L_1$ covers $L_p$ by the map $z\mapsto z^p$ in each fiber, as can be seen by the \v{C}ech construction of $c_1(L_p)$. The equivalence of line bundles over symplectic manifolds and contact manifolds with closed Reeb orbits is given by the Boothby--Wang fibration theorem
\cite[Theorem~3.9]{Blair}.

After discussing two key propositions, we give the proof of Theorem~\ref{thm:main} at the end of this subsection.

We recall the approach of \cite{MRT4}. For any set $X$, the following sets are in bijection:
\[
\maps\bigl(S^1\times X,X\bigr)\leftrightarrow \maps\bigl(S^1, \maps(X,X)\bigr)\leftrightarrow \maps\bigl(X, \maps\bigl(S^1,X\bigr)\bigr).
\]
In particular, let $a\colon S^1\times \omp\to\omp$ be the isometric $S^1$ action of rotation in the fibers of $\omp$. This gives
$a^L\colon \omp\to L\omp$, defined by $a^L(m)(\theta) = a(\theta,m)$, and
$a^I\colon S^1\to\ism(M,g)$, defined by $a^I(\theta)(m) = a(\theta,m)$.
Similarly, for a homotopy,
$F\colon [0,1]\times S^1\times \bm\to \bm$, set
$F^L\colon [0,1]\times \bm\to L\bm$ by $ F^L(t,m)(\theta) = F(t,\theta, m)$.

\begin{Proposition}\label{customprop:app} Let $(\bm,g)$ be a closed $(2k-1)$-manifold. Then
\begin{gather*}
\int_{\bm} a^{L,*}\ocs^{\rm W}_{2k-1}(g) \neq 0 \Rightarrow 0\neq \bigl[a^I\bigr]\in \pi_1(\isom(\bm,g))
\qquad\! {\rm and}\qquad\! \bigl[a^I\bigr]^n\! \neq \bigl[a^I\bigr]^m\ {\rm for} \ m\neq n.
\end{gather*}
In particular, if \smash{$\int_{\bm} a^{L,*}\ocs{}^{\rm W}_{2k-1}(g) \neq 0$}, then $\bigl[a^I\bigr]$ has infinite order in $\pi_1(\ism(M,g))$.
\end{Proposition}

The following is proved in Appendix~\ref{app:oldB}, Lemma~\ref{lem:2.3}.

\begin{Proposition}\label{prop:temp}
Let $F\colon [0,1]\times S^1\times \bm\to \bm$ be a smooth homotopy through isometries, i.e.,
$F\bigl(x^0,\theta,\cdot\bigr)\in \ism(\bm,g)$ for all $\bigl(x^0,\theta\bigr)$. Then
\[ {\rm d}_{\itm}F^{L,*}\ocs^{\rm W}_{\kk} =0.\]
\end{Proposition}

Here \smash{$\ocs{}^{\rm W}_{\kk}= \ocs{}^{\rm W}_{\kk}(g)$}.
From now on, we denote $\ism(\bm,g)$ by $\ism(\bm)$.

We now prove Proposition~\ref{customprop:app}, assuming Proposition~\ref{prop:temp}. As a first step, we prove
that if $a_0$ and $a_1$ are homotopic through isometries, then
\begin{equation}\label{a:1}\int_{\bm} a_0^{L,*}\ocs^{\rm W}_{\kk} = \int_{\bm} a_1^{L,*}\ocs^{\rm W}_{\kk}.
\end{equation}
This is just Stokes' theorem: for $i_{x^0}\colon \bm\to [0,1]\times \bm$, $i_{x^0}(\bm) = \bigl(x^0,m\bigr)$,
\begin{align*}\int_{\bm} a_0^{L,*}\ocs^{\rm W}_{\kk} - \int_{\bm} a_1^{L,*}\ocs^{\rm W}_{\kk} &=
\int_{\bm} i_0^*F^{L,*}\ocs^{\rm W}_{\kk} - \int_{\bm} i_1^* F^{L,*}\ocs^{\rm W}_{\kk}\\
&= \int_{\itm} {\rm d}_{\itm} F^{L,*}\ocs^{\rm W}_{\kk} =0,
\end{align*}
by Proposition~\ref{prop:temp}.

Now let $a_n$ be the $n^{\rm th}$ iterate of
 $a$, i.e., $a_n(\theta,m) =
a(n\theta,m)$.
We claim that
\[
\int_{\bm}a_n^{L,*}\ocs^{\rm W}_{\kk} = n\int_{\bm} a^{L,*}\ocs^{\rm W}_{\kk}.
\]
By
(\ref{csg}), every term in \smash{$\ocs^{\rm W}_{\kk}$} is of the
form \smash{$\ints\dot\gamma(\theta)^i f(\theta)_i\, {\rm d}\theta$}, where $f$ is a periodic one-form on the
circle. Each loop $\gamma\in
a^L_1(\bm)$ corresponds to the loop $\gamma(n\cdot)\in a^L_n(\bm)$. Therefore, the term
$\ints\dot\gamma(\theta)^i f(\theta)_i\, {\rm d}\theta$ is replaced by
\[\ints \frac{{\rm d}}{{\rm d}\theta}\gamma(n\theta)^i f(n\theta)_i\, {\rm d}\theta
 = n\int_0^{2\pi} \dot\gamma(\theta)^if(\theta)_i\, {\rm d}\theta.\]
Thus \smash{$\int_{\bm}a_n^{L,*}\ocs{}^{\rm W}_{\kk} = n\int_{\bm}a^{L,*}\ocs{}^{\rm W}_{\kk}$}.
According to (\ref{a:1}),
$a_n$ and $a_m$ are not homotopic through iso\-metries. Since $\bigl[a^I\bigr]^n = \bigl[a_n^I\bigr]$,
 the $\bigl[a^I\bigr]^n\in
\pi_1(\ism(\bm))$
are all distinct. This proves Proposition~\ref{customprop:app}.

\begin{Remark}
The tricky point of Proposition~\ref{prop:temp} is determining which class of homotopies $F$ gives \smash{\raisebox{-0,3pt}{$ {\rm d}_{\itm}F^{L,*}\ocs{}^{\rm W}_{\kk} =0$}}, or equivalently
\begin{equation}\label{eq:fdlm}F^{L,*} {\rm d}_{LM}\ocs^{\rm W}_{\kk} =0,
\end{equation}
by Lemma~\ref{lem:2.1}. In Lemma~\ref{lem:2.2}, we give a general local formula for
\smash{$F^{L,*} {\rm d}_{LM}\ocs{}^{\rm W}_{\kk}$}. In~\cite{MRT4}, we incorrectly assumed that (\ref{eq:fdlm}) holds when $F$ is a smooth homotopy through diffeomorphisms. Indeed, the proof of~(\ref{eq:fdlm}) when $F$ is a smooth homotopy through isometries (Lemma~\ref{lem:2.3}) is not valid for diffemorphisms; the proof uses the crucial transformation formula~(\ref{eq:11}) for isometries.
\end{Remark}

By Proposition~\ref{customprop:app}, we want to compute \smash{$a^{L,*}\ocs{}^{\rm W}_{4n+1}$} in
local coordinates.
In our setting,~$a^L(m)$ is the loop $\gamma=\gamma_m$ given by the fiber $M_{p,m}$.
We may assume that
$ \bxi= \partial_{x^1}$ is the first element of a local coordinate frame $\{\partial_{x^1},\partial_{x^2},\dots, \partial_{x^{4n+1}}\} = \{\partial_1,\dots,\partial_{4n+1}\}$, and note that $\dot\gamma^\nu = \xi^\nu$.
Then at $\gamma$,~(\ref{csg}) becomes
\begin{equation}\label{CSWlocal}a^{L,*}\ocs^{\rm W}_{4n+1}
 = \ocs^{\rm W}\bigl(a^L_*\partial_1, \dots, a^L_*\partial_{4n+1}\bigr) {\rm d}x^1\wedge \cdots\wedge {\rm d}x^{4n+1}.
\end{equation}

To write~(\ref{CSWlocal}) in local coordinates, we have
\begin{align}
\bigl(a^L_* \partial_i\bigr)(\theta) &=\frac{{\rm d}}{{\rm d}t}\biggl|_{t=0} a^L\bigl(tx^i\bigr)(\theta) =
\frac{{\rm d}}{{\rm d}t}\biggl|_{t=0} a\bigl(\theta, tx^i\bigr) =
\frac{{\rm d}}{{\rm d}t}\biggl|_{t=0} a^I(\theta)\bigl(tx^i\bigr)\nonumber\\
&= a^I_*(\theta)(\partial_i)
= \frac{\partial a^I(\theta)^{j_i}}{\partial x^i}\partial_{j_i}
:= \frac{\partial a^{I,j_i}}{\partial x^i}\partial_{j_i}.\label{eq:cyan1}
\end{align}
We now want to substitute~(\ref{csg}) and~(\ref{K_tensor}) into~(\ref{CSWlocal}). Since $a^I(\theta)\in \ism(M,g)$ for fixed $\theta$, we~have
\[ g_{ij}(m) = \bigl(a^{I,*}g\bigr)_{ij}(m) =g_{\lambda\mu}(a(\theta,m)) \frac{\partial a^{I, \lambda}}{\partial x^i}
\biggl|_{_{(\theta,m)}}
\frac{\partial a^{I, \mu}}{\partial x^j}\biggl|_{_{(\theta,m)}},\]
with, e.g., $a^{I,\lambda} = \bigl(a^I(\theta)\bigr)^\lambda$.
Therefore,~(\ref{K_tensor}) becomes
\begin{gather}
R_{ijk\ell}(m) = \bigl(a^{I,*}R\bigr)_{ijk\ell}(m) = R_{\lambda\mu\nu\kappa}(a(\theta,m))
\frac{\partial a^{I, \lambda}}{\partial x^i}\frac{\partial a^{I, \mu}}{\partial x^j}
\frac{\partial a^{I, \nu}}{\partial x^k}\frac{\partial a^{I, \kappa}}{\partial x^{\ell}},\nonumber\\
K_{ji_1\dots i_{4n+1}}(m)
=\bigl(a^{I,*}K\bigr)_{j i_1 i_2\dots i_{4n+1}} \nonumber\\
 \hphantom{K_{ji_1\dots i_{4n+1}}(m)}{}
 = K_{\nu \lambda_1\lambda_2\dots\lambda_{4n+1}}(a(\theta,m))
{\frac{\partial a^{I,\nu}}{\partial x^j}}
\frac{\partial a^{I, \lambda_1}}{\partial x^{i_1}}\frac{\partial a^{I, \lambda_2}}{\partial x^{i_2}} \cdots
\frac{\partial a^{I, \lambda_{4n+1}}}{\partial x^{i_{4n+1}}},\label{eq:11}
\end{gather}
with all partial derivatives evaluated at $(\theta,m)$.

The term $\dot\gamma^\nu = \partial a^{I,\nu}/\partial\theta$ also appears in~(\ref{csg}).
This term simplifies because $a$ is an action: $a(\theta + \theta',m) =a(\theta,a(\theta',m))$, where angles are added mod $2\pi$. Then,
\begin{align}
\frac{\partial a^{I,\nu}}{\partial\theta} &= \frac{{\rm d}}{{\rm d}\theta'}\biggl|_{\theta'=0}
a^\nu\bigl(\theta+\theta',m\bigr) = \frac{{\rm d}}{{\rm d}\theta'}\biggl|_{\theta'=0} a^\nu\bigl(\theta,a\bigl(\theta',m\bigr)\bigr)
\nonumber\\
&= \frac{\partial a^{I,\nu}}{\partial x^k}\biggl|_{(\theta,m)}\frac{\partial a^k}{\partial \theta'}\biggl|_{\theta'=0}
= \frac{\partial a^{I,\nu}}{\partial x^k}\xi^k,\label{eq:cyan2}
\end{align}
again evaluated at $(\theta,m)$.

Plugging~(\ref{eq:cyan1})--(\ref{eq:cyan2}) into~(\ref{CSWlocal}) at $a(\theta,m)$,
and using the definition of \smash{$\ocs^{\rm W}$} in (\ref{csg}) and~(\ref{K_tensor}), we obtain at $m$
\begin{align}
a^{L,*}\ocs^{\rm W}_{4n+1}
&= \ocs^{\rm W}\bigl(a^L_*\partial_1, \dots, a^L_*\partial_{4n+1}\bigr) {\rm d}x^1\wedge\cdots \wedge {\rm d}x^{4n+1}\nonumber\\
&= \ocs^{\rm W} (\partial_{j_1},\partial_{j_2},\dots,\partial_{j_{4n+1}} )
 \frac{\partial a^{I,j_1}}{\partial x^1} \frac{\partial a^{I,j_2}}{\partial x^2}\cdots
 \frac{\partial a^{I,j_{4n+1}}}{\partial x^{4n+1}}{\rm d}x^1\wedge\cdots \wedge {\rm d}x^{4n+1}\nonumber\\
 &= \frac{4n+1}{2^{4n-1}}\biggl(\int_{S^1} K_{\nu j_1 j_2 \dots j_{4n+1}}(a(m,\theta))
 \frac{\partial a^{I,\nu}}{\partial \theta}
 \frac{\partial a^{I,j_1}}{\partial x^1}\cdots
 \frac{\partial a^{I,j_{4n+1}}}{\partial x^{4n+1}}\, {\rm d}\theta\biggr)\nonumber\\
 &\quad \times {\rm d}x^1\wedge\cdots \wedge {\rm d}x^{4n+1}\nonumber\\
 &= \frac{4n+1}{2^{4n-1}}\biggl(\int_{S^1} K_{\nu j_1 j_2 \dots j_{4n+1}}(a(m,\theta))
 \frac{\partial a^{I,\nu}}{\partial x^k}\xi^k
 \frac{\partial a^{I,j_1}}{\partial x^1}\cdots
 \frac{\partial a^{I,j_{4n+1}}}{\partial x^{4n+1}}\, {\rm d}\theta\biggr)\nonumber\\
 &\quad \times {\rm d}x^1\wedge\cdots \wedge {\rm d}x^{4n+1}\nonumber\\
 &= \frac{4n+1}{2^{4n-1}}\biggl(\int_{S^1} K_{k 1 2\dots 4n+1}(m)\xi^k(m) \, {\rm d}\theta\biggr){\rm d}x^1\wedge\cdots \wedge {\rm d}x^{4n+1}\nonumber\\
 &= \frac{2\pi(4n+1)}{2^{4n-1}} K_{k 1 2\dots 4n+1}(m) \xi^k(m)\, {\rm d}x^1\wedge\cdots \wedge {\rm d}x^{4n+1}.\label{eq:cyan3}
\end{align}
Since $\xi = \partial_1$, we have $\xi^k = \delta^{1k}$, so (\ref{eq:cyan3}) becomes
\begin{equation}\label{eq:cyan4}
a^{L,*}\ocs^{\rm W}_{4n+1}(m) = \frac{2\pi(4n+1)}{2^{4n-1}} K_{1, 1 2\dots 4n+1}(m)\, {\rm d}x^1\wedge\cdots \wedge {\rm d}x^{4n+1}.
\end{equation}

In the definition of $K_{1,12\dots 4n+1}$ in (\ref{K_tensor}) in terms of the curvature tensor, we can substitute~the formulas in Lemma~\ref{lem1} to write the top form
\smash{$a^{L,*}\ocs{}^{\rm W}_{4n+1}$} on $\omp$ as a polynomial in $p$ with curvature expressions as coefficients:
\begin{equation}\label{eq:11aa}
a^{L,*}\ocs^{\rm W}_{4n+1}(m) = \frac{2\pi(4n+1)}{2^{4n-1}}\sum_{q = 1}^{4n+2} S_{4n+1,q}(m) p^q {\rm d}x^1\wedge\dots\wedge {\rm d}x^{4n+1}.
\end{equation}

We focus on the top term.

\begin{Claim}\label{claim1}
For fixed $m$,
\begin{align}
S_{4n+1,4n+2} &= S_{4n+1,4n+2}(m) \nonumber\\
&=(-1)^{n+1} 2^{2n+1}(2n+1)
\sum_{\sigma \in \mathfrak S_{4n}} \sgn(\sigma) [J_{\sigma_1 \sigma_2} \cdots J_{\sigma_{4n-1} \sigma_{4n}} ](m),\label{Seqn}
\end{align}
where $\mathfrak S_{4n}$ is the permutation group of $\{1,\dots, 4n\}$.
\end{Claim}

This is proved in Appendix~\ref{app:claim}.

\begin{proof}[Proof of Theorem~\ref{thm:main}, assuming Proposition~\ref{prop:temp} and Claim~\ref{claim1}.]
By Lem\-ma~\ref{lem:Y.one}, the \linebreak right-hand side of (\ref{Seqn}) is a nonzero multiple of the symplectic volume form.
As a result, (\ref{eq:11aa}) and Claim~\ref{claim1} imply
\[0\neq \int_{\omp}a^{L,*}\ocs^{\rm W}_{4n+1}\]
for $|p| \gg 0$.
By Proposition~\ref{customprop:app}, $\bigl[a^I\bigr]$ has infinite order in $\pi_1(\isom(\omp))$ for $|p| \gg 0$.
Thus Theorem~\ref{thm:main} follows.
\end{proof}

Theorem~\ref{thm:main} applies to symplectic manifolds $\om$ of dimension $4n$. The crucial Claim~\ref{claim1} does not hold if ${\rm dim}\bigl(\om\bigr) = 4n+2$. In Appendix~\ref{appb}, we prove that $S_{7,8} = 0$ for ${\rm dim}(M) = 6$, and the proof extends directly to $S_{4n+3, 4n+4}$.

\subsection{A computer verification}\label{section3.2}

Using the code at {\tt egison.org}, we obtain the following results for
$S_{4n+1, 4n+2}$
in Theorem~\ref{thm:main}~\cite{Egi1}.

\begin{center}\begin{tabular}{|c||c|c|c|}
\hline
dim(M) & 4&\ 6\ \ &8\\
\hline
$S_{4n+1,4n+2}$ & $-192$&0 &61440\\
\hline
\end{tabular}
\end{center}

In this pointwise calculation, we have
put the almost complex structure into the normal form
\[J = \begin{pmatrix} \hphantom{-} 0& \rm I\\ -\rm I&0\end{pmatrix}.\]
This agrees with (\ref{3.5add12}): for ${\rm dim}(M)=4$ (i.e., $n=1$), $\sum _{\sigma_0=0} J_{\sigma_1\sigma_2}J_{\sigma_3\sigma_4} =
-2!2^2$, so $S_{4n+1,4n+2} = (-1)^2 2^3(3) (-8) =-192$; for ${\rm dim}(M)=8$ (i.e., $n=2$), the corresponding
sum over permutations gives $4!2^4$, so
$S_{4n+1,4n+2} = (-1)^3 2^5 5 (384)= 61440$.

\section{The Kodaira--Thurston example}\label{section4}

We calculate explicitly the WCS class for the example independently given by Kodaira \cite{Kodaira} and
Thurston \cite{Thurston} of a non-K\"ahler symplectic manifold \smash{$\om^4$}. By putting an explicit Riemannian metric $g_p$ on $M_p$, we can compute that
$\pi_1({\rm Isom}(\omp, g_p))$ are infinite for all $p\in \Z$.

\subsection{The metric}
The orbit space $\om$ is
a $T^2$ fibration over $T^2$. To construct $\om$, we take coordinates $\theta_1, \theta_2,\theta_3,\theta_4\in [0,1]$. The base~$T^2 $ has coordinates $\theta_1$, $\theta_2$, where we glue $\theta_1$, $\theta_2$ as usual to get
a torus. For the fiber $T^2$, we take the linear transformation
$\left(\begin{smallmatrix} 1&1\\0&1\end{smallmatrix}\right)$
on $\R^2 = \{(\theta_3,\theta_4)\}$ (so now these coordinates are real numbers), which glues
the unit $(\theta_3,\theta_4)$-square to the parallelogram with sides given by the vectors $\vec\theta_3$,
 $\vec \theta_3 + \vec\theta_4$. We do this gluing in the~$\theta_2$ direction, so that $\om$ is given by
 $[0,1]^4$ with the relations/gluings
 \[(0,\theta_2,\theta_3,\theta_4) \sim (1,\theta_2,\theta_3,\theta_4), \qquad
(\theta_1,0,\theta_3,\theta_4) \sim (\theta_1, 1, \theta_3, \theta_3+\theta_4).\]

We claim that the metric
\begin{equation*}
{\rm d}\theta_1^2 + {\rm d}\theta_2^2 + {\rm d}\theta_3^2 -\theta_2{\rm d}\theta_3 {\rm d}\theta_4 + (1+\theta_2){\rm d}\theta_4^2
\end{equation*}
is well-defined on $\om$.
Since $\partial_{\theta_4}$ at $\theta_2=0$ is glued to
$\partial_{\theta_3} + \partial_{\theta_4}$ at $\theta_2 = 1$, this means we must~have%
\begin{align*}
&\langle \partial_{\theta_i},\partial_{\theta_j}\rangle_{(0,\theta_2,\theta_3,\theta_4)}
= \langle \partial_{\theta_i},\partial_{\theta_j}\rangle_{(1,\theta_2,\theta_3,\theta_4)},
\qquad i,j = 1,2,3,4,\\
&\langle \partial_{\theta_i},\partial_{\theta_j}\rangle_{(\theta_1,0,\theta_3,\theta_4)}
= \langle \partial_{\theta_i},\partial_{\theta_j}\rangle_{(\theta_1,1,\theta_3,\theta_4)},\qquad i, j = 1, 2, 3,\\
&\langle \partial_{\theta_i},\partial_{\theta_4}\rangle_{(\theta_1,0,\theta_3,\theta_4)}
= \langle \partial_{\theta_i},\partial_{\theta_3} + \partial_{\theta_4}\rangle_{(\theta_1,1,\theta_3,\theta_4)},\qquad
i = 1,2,3,\ j = 4,\\
&\langle \partial_{\theta_4},\partial_{\theta_4}\rangle_{(\theta_1,0,\theta_3,\theta_4)}
= \langle \partial_{\theta_3} + \partial_{\theta_4},\partial_{\theta_3} + \partial_{\theta_4}\rangle_{(\theta_1,1,\theta_3,\theta_4)}.
\end{align*}
Since the metric is independent of $\theta_1\in [0,1]$, the first equation holds; since
the metric is independent of $\theta_2$ for $i,j = 1,2,3$, the second equation holds.
 For the third equation, the left-hand side is $0$; the right-hand side is also $0$ for $i=1,2$, and for $i=3$ we get
\[\langle \partial_{\theta_i},\partial_{\theta_3} + \partial_{\theta_4}\rangle_{(\theta_1,1,\theta_3,\theta_4)}
= 1 - (\theta_2 = 1) = 0.\]
For the last equation, the left-hand side is $1$, and the right-hand side is
\[\langle\partial_{\theta_3},\partial_{\theta_3}\rangle_{\theta_2=1} +
2\langle\partial_{\theta_3},\partial_{\theta_4}\rangle_{\theta_2=1} +
\langle\partial_{\theta_4},\partial_{\theta_4}\rangle_{\theta_2=1}
= 1 +2 (-1) + 2 = 1.\]
(Since $g_{33} = 1$ is independent of $\theta_2$, from the gluing
$\partial_{\theta_3}|_{\theta_2=0} = \partial_{\theta_3}|_{\theta_2=1}$, $\partial_{\theta_4}|_{\theta_2=0} = (\partial_{\theta_3}+\partial_{\theta_4})|_{\theta_2=1}$,
$g_{34}(\theta_2)$ must satisfy $g_{34}(0) = 0$, $g_{34}(1) = -1$ and $g_{44}(\theta_2)$ must satisfy $g_{44}(0) = 1$, $g_{44}(1) = 2$. Thus our choice of metric is the simplest one possible.)

As a check, we note that the volume form is
\[ \bigl(1+\theta_2-\theta_2^2\bigr){\rm d}\theta_1\wedge {\rm d}\theta_2 \wedge {\rm d}\theta_3\wedge {\rm d}\theta_4,\]
which is equal at $\theta_2=0$ and $\theta_2 = 1$. It is also positive definite, since $1+\theta_2-\theta_2^2$ has no roots in~$[0,1]$.

\subsection{The compatible AC structure and the new metric}

Given a symplectic form $\omega$ and a Riemannian metric $g$, we want to find an AC structure $J$ and a new metric $\tg$ with the compatibility condition
\begin{equation}\label{oneaa}\omega(u,v) = \tg(Ju,v).
\end{equation}
 The usual procedure is to write
$\omega(u,v) = g(Au, v)$
for some skew-adjoint transformation $A$. (The matrix of $A$ is not necessarily skew-symmetric in the
basis $\{\partial_{\theta_i}\}$, since this basis is only orthogonal at $\theta_2=0$.)
For $A^*$
the adjoint of $A$ with respect to $g$, we set
\begin{equation*}
J = \sqrt{AA^*}^{-1}A = \sqrt{-A^2}^{-1}A,\qquad \tg(u,v) = g\bigl(\sqrt{AA^*}u,v\bigr).
\end{equation*}
It is easy to check that $J^2 = -1$ and that (\ref{oneaa}) holds. Note that
\[
\tg(u,v) = g\bigl((AA^*)^{1/4}u,(AA^*)^{1/4}v\bigr)
\]
is positive definite.

We take the symplectic form $\omega = {\rm d}\theta_1\wedge {\rm d}\theta_2 + \kappa {\rm d}\theta_3\wedge {\rm d}\theta_4$, $\kappa\in \Z\setminus \{0\}$, so $\bigl(\om,\omega\bigr)$ is integral. \big(For $\kappa <0$, $\omega^2$ is the volume form for the reverse of the standard orientation.\big) For the metric~$g$, we first have to compute $A$. The compatibility condition (\ref{oneaa}) is equivalent to
\begin{equation*}
\omega_{ij} = A^k_i g_{kj}.
\end{equation*}
A straightforward calculation gives
\[A = \begin{pmatrix} 0&1&0&0\\
-1&0&0&0\\
0&0& \dfrac{\theta_2\kappa}{1+\theta_2-\theta_2^2}& \dfrac{\kappa}{1+\theta_2-\theta_2^2}\vspace{1mm}\\
0&0&\dfrac{(-1-\theta_2)\kappa}{1+\theta_2-\theta_2^2}& \dfrac{-\theta_2\kappa}{1+\theta_2-\theta_2^2}
\end{pmatrix}\]

We now have to compute $\sqrt{AA^*}$. From (\ref{oneaa}) and
\[\omega(u,v) = -\omega(v,u) = -g(Av,u) = g(-A^*u,v),\]
we get $A^* = -A$. Thus
\begin{gather*}
AA^* =
\begin{pmatrix} 1&0&0&0\\ 0&1&0&0\\
0&0& \dfrac{\kappa^2}{1+\theta_2-\theta_2^2}&0\\
0&0&0&\dfrac{\kappa^2}{1+\theta_2-\theta_2^2}\\
\end{pmatrix}\\
\qquad
\Longrightarrow \sqrt{AA^*} =
 \begin{pmatrix}1&0&0&0\\ 0&1&0&0\\
0&0& \dfrac{\kappa}{\bigl(1+\theta_2-\theta_2^2\bigr)^{1/2}}&0\\
0&0&0&\dfrac{\kappa}{\bigl(1+\theta_2-\theta_2^2\bigr)^{1/2}}\\
\end{pmatrix},
\end{gather*}
and{\samepage
\begin{equation*}
J = \sqrt{AA^*}^{-1}A =
\begin{pmatrix} 0&1&0&0\\
-1&0&0&0\\
0&0& \dfrac{\theta_2}{\bigl(1+\theta_2-\theta_2^2\bigr)^{1/2}}& \dfrac{1}{\bigl(1+\theta_2-\theta_2^2\bigr)^{1/2}}\vspace{1mm}\\
0&0&\dfrac{-1-\theta_2}{\bigl(1+\theta_2-\theta_2^2\bigl)^{1/2}}&\dfrac{-\theta_2}{\bigl(1+\theta_2-\theta_2^2\bigr)^{1/2}}
\end{pmatrix}.
\end{equation*}
Note that $J$ is independent of $\kappa$.}

To compute $\tg$, we have
\begin{align*}
\sqrt{AA^*}\begin{pmatrix} u_1\\u_2\\u_3\\u_4\end{pmatrix} &= \begin{pmatrix} u_1\\u_2\\
\dfrac{\kappa}{\bigl(1+\theta_2-\theta_2^2\bigr)^{1/2}}u_3 \\
\dfrac{\kappa}{\bigl(1+\theta_2-\theta_2^2\bigr)^{1/2}}u_4
\end{pmatrix} \\
&\Longrightarrow \tg = g\bigl(\sqrt{AA^*}\cdot,\cdot\bigr) = \begin{pmatrix} 1&0&0&0\\
0&1&0&0\\
0&0& \dfrac{\kappa}{\bigl(1+\theta_2-\theta_2^2\bigr)^{1/2}}& \dfrac{-\theta_2\kappa}{\bigl(1+\theta_2-\theta_2^2\bigr)^{1/2}}\\[2mm]
0&0&\dfrac{-\theta_2\kappa}{\bigl(1+\theta_2-\theta_2^2\bigr)^{1/2}}& \dfrac{(1+\theta_2)\kappa}{\bigl(1+\theta_2-\theta_2^2\bigr)^{1/2}}
\end{pmatrix}.
\end{align*}

We now use $\tg$ to define $g = g_p$ on $M_p$ as in (\ref{eq:gp}).

\subsection{The top WCS form}

Let $\{e_0,\dots,e_4\}$ be a local orthonormal frame of $M_p$, with $e_0 = \bxi$.
By (\ref{csg}) with $k=3$,
\begin{align}
a^{L,*}\widetilde{{\rm CS}}^{\rm W}_5(e_0,\dots,e_4) &= \widetilde{{\rm CS}}^{\rm W}_5 \bigl(a^L_*e_0,\dots, a^L_* e_4\bigr)\nonumber\\
&=\frac{3}{2}\int_{S^1}
\sum_{\sigma \in {\mathfrak S_5}}
\sgn(\sigma) \barcurvature{\sigma_0 \ell_1 0}{r}
\barcurvature{\sigma_1 \sigma_2 \ell_2}{\ell_1}
\barcurvature{\sigma_3\sigma_4 r}{\ell_2}{\rm d}\theta_0,\label{CSThurston}
\end{align}
where $\mathfrak S_5$ is the permutation group on $\{0,1,2,3,4\}$, $\ell_1, \ell_2, r\in \{0,1,2,3,4\}$,
$\sigma=(\sigma_0, \sigma_1, \sigma_2,\allowbreak \sigma_3,\sigma_4)$, and $\theta_0$ is the fiber coordinate with $\partial_{\theta_0} = \bxi$. We have used that $a^L(m)$ is the circle fiber of $m\in \omp$, so $\dot\gamma$ in (\ref{csg}) equals $\bxi$. Thus
the integral over $S^1$ is the integral over the circle fiber in $\omp$.

Set $\beta = \beta(\theta_2) = 1+ \theta_2-\theta_2^2$.

\begin{Proposition}\label{Thurston theorem}
We have
\begin{equation}\label{thm:thurston}
\int_{\omp}
a^{L,*}\widetilde{{\rm CS}}^{\rm W}_5=
\frac{3\kappa\pi^2 p^{2}}{8}
\int_0^1\bigl(3072p^4 -640p^2\beta^{-2} -25\beta^{-4}\bigr){\rm d}\theta_2.
\end{equation}
\end{Proposition}

\begin{proof} We explain the constants on the right-hand side of (\ref{thm:thurston}). By the construction of~$\bg$,~$a$~acts via isometries on $\omp$. As in (\ref{eq:cyan3}), this makes the integrand in (\ref{CSThurston}) independent of $\theta_0$, so the integral
is replaced with a factor of $2\pi$. Thus
\begin{align*} \int_{\omp} a^{L,*}\widetilde{{\rm CS}}^{\rm W}_5 &=
 \int_{\omp} a^{L,*}\widetilde{{\rm CS}}^{\rm W}_5(e_0,\dots,e_4) e^0\wedge\cdots\wedge e^4\\
 &=\frac{2\pi\cdot 3}{2} \int_{\omp}
 \sum_{\sigma \in {\mathfrak S_5}}
\sgn(\sigma) \barcurvature{\sigma_0 \ell_1 0}{r}
\barcurvature{\sigma_1 \sigma_2 \ell_2}{\ell_1}
\barcurvature{\sigma_3\sigma_4 r}{\ell_2}\dvol.
\end{align*}
We now switch to the coordinates $\{\theta_0,\dots,\theta_4\}$, so $\br$ is now computed in these coordinates, and $\dvol = \kappa\, {\rm d}\theta_0\wedge\dots
\wedge {\rm d}\theta_4$.
The integrand is again independent of the point in the fiber, so the integral over the fiber just detects the length of the fiber. By the construction of $\bg$, the fiber in~$M_1$ has length \smash{$2\pi
=\int_0^{2\pi} |\bxi|$}.
Since
$g = g_p$ involves $\be = \be_p$, and since $\bxi = \bxi_p$ has $\be_p(\bxi_p) = 1$,
each fiber has $g_p$ length \smash{$\int_0^{2\pi} |\xi_p|_{g_p}{\rm d}\theta = 2\pi$} independent of $p$. Therefore,
\begin{align*} \int_{\omp} a^{L,*}\widetilde{{\rm CS}}^{\rm W}_5 &=
3\pi\cdot 2\pi\kappa
\int_M \sum_{\sigma \in {\mathfrak S_5}}
\sgn(\sigma) \barcurvature{\sigma_0 \ell_1 0}{r}
\barcurvature{\sigma_1 \sigma_2 \ell_2}{\ell_1}
\barcurvature{\sigma_3\sigma_4 r}{\ell_2}{\rm d}\theta_1\wedge\cdots \wedge {\rm d}\theta_4.
\end{align*}
Thus the proposition
follows if
\begin{align}
&\int_M \sum_{\sigma \in {\mathfrak S_5}}
\sgn(\sigma) \barcurvature{\sigma_0 \ell_1 0}{r}
\barcurvature{\sigma_1 \sigma_2 \ell_2}{\ell_1}
\barcurvature{\sigma_3\sigma_4 r}{\ell_2}{\rm d}\theta_1\wedge\cdots \wedge {\rm d}\theta_4 \nonumber \\
&\qquad= \frac{p^2}{16} \int_0^1\bigl(3072p^4 -640p^2\beta^{-2} -25\beta^{-4}\bigr){\rm d}\theta_2.\label{tobeshown}
\end{align}
The long calculation of (\ref{tobeshown}) is in \cite{Egi3}. This result is verified by the computer calculations in a~file
at~\cite{Egi1}.
\end{proof}

Since the top coefficient of $p$ is nonzero, we conclude from Theorem~\ref{thm:main} that
$\pi_1({\rm Isom}(\omp,g_p))$ is infinite for $|p| \gg 0$. We will improve this to all $p$ as follows.

\begin{Theorem}
$\pi_1({\rm Isom}(\omp,g_p))$ is infinite for all $p$.
\end{Theorem}

\begin{proof} For $p=0$, this follows from $\omz = \om \times S^1$ (cf.\ \cite[Remark~3.2]{MRT4}).
For $p\neq 0$, by (\ref{thm:thurston}) and~(\ref{tobeshown}), it suffices to show that
\begin{equation}\label{tbs}
 \int_0^1\bigl(3072p^4 -640p^2\beta^{-2} -25\beta^{-4}\bigr){\rm d}\theta_2\neq 0
\end{equation}
for $p\in \Z$. Either by a direct calculation or by Wolfram Alpha, we get (for $\theta= \theta_2$)
\begin{align*}
&\int \beta^{-2}{\rm d}\theta= \frac{2\theta-1}{5(1+\theta-\theta^2)} -
\frac{2(\ln\bigl( -2\theta+ \sqrt{5} + 1\bigr)}{5\sqrt{5}} + \frac{2\bigl(\ln( 2\theta+ \sqrt{5} - 1\bigr)}{5\sqrt{5}} + C, \\
&\int \beta^{-4} {\rm d}\theta=-\frac{1}{375}\biggl( \frac{1}{(1+\theta-\theta^2)^3}\bigl(-60\theta^5+150 \theta^4
+50x^3-225x^2-75\theta+ 80\bigr) \\
&\hphantom{\int \beta^{-4} {\rm d}\theta=-\frac{1}{375}\biggl(}{}
+12\sqrt{5}\bigl(\ln\bigl(2\theta+ \sqrt{5}- 1\bigr)-\ln\bigl(-2\theta+\sqrt{5}+1\bigr)\bigr)
\biggr) +C.
\end{align*}
The definite integrals are
\begin{align*}
&\int_0^1 \beta^{-2}{\rm d}\theta= \frac{2}{25}\bigl( 5+4\sqrt{5}\coth^{-1}\bigl(\sqrt{5}\bigr)\bigr),\\
& \int_0^1 \beta^{-4}{\rm d}\theta=\frac{16}{375}\bigl( 10+3\sqrt{5}\coth^{-1}\bigl(\sqrt{5}\bigr)\bigr).
 \end{align*}
 Plugging this into (\ref{tbs}), we must show that
\[10\bigl(-1 - 24p^2 + 288 p^4\bigr) -3\sqrt{5}\bigl(1+64p^2\bigr) \coth^{-1}\bigl(\sqrt{5}\bigr)\neq 0.\]
This quadratic equation in $p^2$ has solutions $p \approx \pm 0.159514 {\rm i}, \pm 0.424868$. Since there are no integral solutions, the theorem
 follows.
 \end{proof}

 A second computer program verifying these calculations is in \cite{Egi1}.

\section{The K\"ahler case}\label{section5}

In this section, we prove that the lowest order term in the WCS form has
a geometric/topological interpretation on
K\"ahler manifolds (Proposition~\ref{app1}); this appears to fail for general symplectic manifolds. We use this result to give non-vanishing results for a type of Wodzicki--Pontryagin form in dimension $4k+2$ on loop spaces (Theorem~\ref{nonv}).
This is an infinite-dimensional phenomenon, as the finite-dimensional version of these forms vanish.
As noted in (\ref{dwcs}), the Wodzicki--Pontryagin classes vanish in
$H^{4k+2}(LM)$; this non-vanishing of the representative forms gives the first known examples in arbitrarily high dimensions.

We start with a result about the real cohomology of loop spaces. This is not used later, but we think it is of general interest.

For a manifold $N$, define the ring homomorphism $L\colon \Lambda^*(N)\to \Lambda^*(LN)$, $ \delta\mapsto \delta_L$, by
\[\delta_L(X_1,\dots,X_{k}) = \delta(X_1(0),\dots, X_{k}(0)).\]
Let $a\colon S^1\times N\to N$ be an $S^1$ action; $a$ can be the trivial action, so there is no restriction on~$N$. Then
\begin{equation}\label{ident}a^{L,*}\circ L = {\rm Id}.
\end{equation}
To see this, take $v\in T_pN$ and a curve $\gamma(s)$ tangent to $v$ at $p$. Then $a^L_{*}(v)
= ({\rm d}/{\rm d}s)|_{s=0} a^L(\gamma(s))$, a vector field along the loop \smash{$a^L(p)(\theta)$}. Since \smash{$a^L(\gamma(s)(0) = \gamma(s)$}, we get \smash{$a^L_{*}(v)(\theta=0) = v$}. Thus
\begin{align*}a^{L,*}\delta_L(X_1,\dots, X_{r}) &= \delta_L\bigl(a^L_{*}X_1,\dots,a^L_{*}(X_{r})\bigr) = \delta
\bigl(\bigl[a^L_{*}X_1\bigr](0),\dots,\bigl[a^L_{*}(X_{r})\bigr](0)\bigr)\\
&= \delta(X_1,\dots,X_r).
\end{align*}
By Lemma~\ref{lem:2.1}, $a^{L,*}\colon \Lambda^*(LN) \to \Lambda^*(N)$ induces a map $\bigl(a^{L,*}\bigr)^*$, just denoted
$a^{L,*}\colon H^*(LN,\R)\to H^*(N,\R)$, on de Rham cohomology.

\begin{Lemma}
The ring homomorphism $L$ induces an injection $L^*\colon H^k(N,\R)\hookrightarrow H^k(LN,\R)$ for all $k$, and
$a^{L,*}\colon H^k(LN,\R)\to H^k(N,\R)$ is surjective for all $k$.
\end{Lemma}

\begin{proof} We have to check that $L$ induces a map on cohomology. As in Lemma~\ref{lem:2.1}, we have
\begin{gather*}\lefteqn{[ ( {\rm d}_{LN}\circ L) \delta](X_1,\dots, X_{k+1})}\\
\qquad{} = \sum_i (-1)^{i-1}X_i\bigl(\delta_L\bigl(X_1,\dots,\hat X_i,\dots X_{k+1} \bigr)\bigr)\\
 \qquad\quad{} + \sum_{i< j}(-1)^{i+j} \delta_L\bigl([X_i,X_j], X_1,\dots, \hat X_i,\dots, \hat X^j,\dots, X_{k+1}\bigr) \\
\qquad{}= \sum_i (-1)^{i-1}X_i\bigl(\delta\bigl(X_1(0),\dots,\hat X_i(0),\dots X_{k+1}(0) \bigr)\bigr)\\
\qquad\quad{} + \sum_{i< j}(-1)^{i+j} \delta\bigl([X_i,X_j](0), X_1(0)\dots, \hat X_i(0),\dots, \hat X^j(0),\dots, X_{k+1}(0)\bigr),\\
[(L\circ {\rm d}_N)\delta](X_1,\dots, X_{k+1})\\
\qquad{}= {\rm d}_N\delta(X_1(0),\dots, X_{k+1}(0)) \\
\qquad{}= \sum_i (-1)^{i-1}X_i(0)\bigl(\delta(X_1(0),\dots,\hat X_i(0),\dots X_{k+1}(0) )\bigr)\\
\qquad\quad{} + \sum_{i< j}(-1)^{i+j} \delta\bigl([X_i,X_j](0), X_1(0)\dots, \hat X_i(0),\dots, \hat X^j(0),\dots, X_{k+1}(0)\bigr).
\end{gather*}
Let $\gamma_s(t)$ be a family of loops with tangent vector $X_i\in T_{\gamma_0}LN$. Extend the $X_j$ to vector fields near $\gamma = \gamma_0$. Then
\begin{align*}X_i \bigl( \delta\bigl(X_1(0),\dots,\hat X_i(0),\dots X_{k+1}(0) \bigr)\bigr) &= \left(\frac{{\rm d}}{{\rm d}s}\biggr|_{s=0} \right) \delta_{\gamma(s)(0)}
\bigl(X_1(0),\dots,\hat X_i(0),\dots X_{k+1}(0)\bigr)\\
& = X_i(0) \bigl( \delta\bigl(X_1(0),\dots,\hat X_i(0),\dots X_{k+1}(0) \bigr)\bigr).
\end{align*}
It follows that $ {\rm d}_{LN}\circ L = L\circ {\rm d}_N$, so
$L\colon \Lambda^*(N)\to\Lambda^*(LN)$ induces $L^*\colon H^*(N) \to H^*(LN)$. Then~(\ref{ident}) implies
$ a^{L,*} L^*= {\rm Id}$, which gives the injectivity of $L^*$ and the surjectivity of $a^{L,*}$.
\end{proof}

In contrast to this general cohomological
result, our goal is to obtain
information on the WCS forms on $M_p$ from the Pontryagin forms on $M$.

Let $\bigl(\om,\omega\bigr)$ be an integral K\"ahler manifold of real dimension $4n$.
The K\"ahler form $\omega\in H_2\bigl(\om,\Z\bigr)$ determines the Riemannian metric $\bgg(X,Y) = -\omega(JX,Y)$, where $J$ is the complex structure. The key feature of the K\"ahler case for us is that $\bbn J=0$. Thus in Lemma~\ref{lem1}, the terms~(ii) and~(iv) vanish.

By (\ref{csg}) and Lemma~\ref{lem1}, the WCS forms \smash{$\ocs{}^{\rm W}_{2k-1}$} on~$L\omp$ and their pullbacks
 \smash{$a^{L,*}\ocs{}^{\rm W}_{2k-1}$}
 to~$\omp$ are polynomials in $p^2$:
\begin{align*}
&\ocs^{\rm W}_{2k-1}= \ocs^{\rm W}_{2k-1}(\omp, \bg_p) =\sum_{i=1}^{k}\ocs^{\rm W}_{2k-1, 2i}\ p^{2i}
\in \Lambda^{2k-1}(L\omp),\\
&a^{L,*}\ocs^{\rm W}_{2k-1}= a^{L,*}\ocs^{\rm W}_{2k-1}(\omp, \bg_p) =\sum_{i=1}^{k}a^{L,*}\ocs^{\rm W}_{2k-1, 2i}\ p^{2i}
\in \Lambda^{2k-1}(\omp).\nonumber
\end{align*}
The forms \smash{$\ocs{}^{\rm W}_{2k-1, 2i}$} are curvature expressions independent of $p$. Indeed,
\smash{$\ocs{}^{\rm W}_{2k-1, 2i}$} involves an integration of an $S^1$-invariant, $p$-independent curvature expression over the $S^1$-fiber in $\omp$. Each fiber has length $2\pi$ independent of $p$ (see the proof of Theorem~\ref{Thurston theorem}),
so \smash{$\ocs{}^{\rm W}_{2k-1, 2i}$} is independent of $p$.

In the proof of the main Theorem~\ref{thm:main}, we computed the highest power of $p$ in~\smash{$a^{L,*}\ocs{}^{\rm W}_{4n+1}$} in~terms of the symplectic structure (see (\ref{Seqn})). In contrast, the lowest power of $p$ in \smash{$a^{L,*}\ocs{}^{\rm W}_{4k+1}$}, for any $k$, contains Pontryagin-type form information.

\begin{Proposition} \label{app1}Let $\pi\colon \omp\to \om$ be the fibration. For $\bxi$ the unit tangent vector to the fibers of
$\pi$,
\begin{align*}\imath_{\bxi} a^{L,*}\ocs^{\rm W}_{4k+1,2} &=
(2k+1)2\cdot \pi^*\tr\bigl(\bo_{\om}^{2k}\bigr)
=
(-1)^k (4k+2)(2\pi)^{2k+1}(2k)!
\cdot \pi^*\tilde p_k\bigl(\bo_{\om}\bigr).
\end{align*}
\end{Proposition}

Here \smash{$\bo = \bo_{\om}$} is the curvature of the K\"ahler metric, and $i_\xi$ is interior product with $\xi$.
 The proof is in Appendix~\ref{app:E}.

There are no finite-dimensional Pontryagin forms in dimensions $4k+2$, because $\Omega^{2k+1}$ is skew-symmetric and hence
$\operatorname{Tr}\bigl(\Omega^{2k+1}\bigr) =0$ in finite dimensions. For the $s=1$ Sobolev connection on~$L\omp$, the curvature is a skew-adjoint zeroth order $\Psi$DO. The top order symbol is easily seen to be skew-symmetric, but the minus one order symbol used to compute the Wodzicki residue need not be skew-symmetric. Thus we can define Wodzicki--Pontryagin forms in dimensions~${4k+2}$.

\begin{Definition} The Wodzicki--Pontryagin form \smash{$\tp^{\rm W}_{k+\frac{1}{2}}\in \Lambda^{4k+2}(L\omp)$} is
\[\tp^{\rm W}_{k+\frac{1}{2}} :=
\resw\bigl(\Omega^{2k+1}\bigr).\]
\end{Definition}

 We now give examples where these Wodzicki--Pontryagin forms are nonzero. The proof depends on the following observations:
 \begin{enumerate}\itemsep=0pt
 \item By (\ref{dwcs}),
the Wodzicki--Pontryagin classes satisfy
 \begin{equation*}
 \tp^{\rm W}_{\frac{k}{2}} = -{\rm d}^{LM_p}\ocs^{\rm W}_{2k-1}\in \Lambda^{2k}(LM_p).
 \end{equation*}
 \item If $\tp^{\rm W}_{\frac{k}{2}} \equiv 0$, then,
 \[F^{L,*} {\rm d}^{LM_p}\ocs^{\rm W}_{2k-1} = {\rm d}_{\itm}F^{L,*}\ocs^{\rm W}_{2k-1} =0,\]
for any smooth map
$F\colon [0,1]\times S^1\times \bm\to \bm$, with
$F^L\colon [0,1]\times \bm\to L\bm$ defined by $ F^L(t,m)(\theta) = F(t,\theta, m)$. In particular, $F$ could be a homotopy through isometries, conformal diffeomorphisms, smooth simple homotopy equivalences, or any subgroup $\mathcal{G}$ of the diffeomorphism group $\diff(M)$.
\item Thus we can apply
Proposition~\ref{prop:temp} in the proof of Proposition~\ref{customprop:app} and conclude that \smash{$a^{L,*}\ocs{}^{\rm W}_{2k-1}(g) \neq 0$} implies $\pi_1(\mathcal{G})$ is infinite, a suspiciously strong result.
 \end{enumerate}

For the example, let $\calL_p = S^5/\Z_p$ be the lens space given by the identification $z \sim {\rm e}^{2\pi {\rm i}/p}z$ for $
z\in S^5$.
By \cite[Proposition 3.14]{MRT4}, $\calL_p$ is diffeomorphic to \smash{$M_p := \overline{\mathbb{CP}^2}_p$}, where the base space $\mathbb{CP}^2$ is of course K\"ahler.

\begin{Theorem}\label{nonv}
\smash{$\tp^{\rm W}_{\frac{3}{2}} \in \Lambda^{6}L(\calL_p)$} is not identically zero.
\end{Theorem}

\begin{proof}
As above, \smash{$\tp^{\rm W}_{\frac{3}{2}} = -{\rm d}^{LM_p}\ocs^{\rm W}_{5}\in \Lambda^6(L\calL_p)$}.
If \smash{$\tp^{\rm W}_{\frac{3}{2}}\equiv 0$}, then
\[F^{L,*} {\rm d}^{LM_p}\ocs^{\rm W}_{5} = {\rm d}_{\itm}F^{L,*}\ocs^{\rm W}_{5} =0,\]
for
$F\colon [0,1]\times S^1\times \bm\to \bm$
a homotopy through smooth simple homotopy equivalences.
By the calculation in
\cite[equation~(3.8)]{MRT4} \big(with the typo \smash{$\int_{\overline {\mathbb{CP}^2}_1}$} replaced with
\smash{$\int_{\overline {\mathbb{CP}^2}_p}$}\big), we obtain \smash{$\int_{\calL_p} a^{L,*}\ocs{}^{\rm W}_5\neq 0$} for $p>1$. Thus
Proposition~\ref{customprop:app} implies \smash{$\pi_1(C^\infty{\rm Aut}(\calL_p))$} is infinite, where $C^\infty{\rm Aut}(\calL_p)$ is the space of smooth simple homotopy equivalences. Since $C^\infty{\rm Aut}(\calL_p)$ is an open subset of $C^\infty(\calL_p,\calL_p)$, and similarly for the corresponding $C^0$ spaces,
the inclusion of~$C^\infty{\rm Aut}(\calL_p)$ into $C^0{\rm Aut}(\calL_p)$ is a homotopy equivalence \cite[Theorem~16]{palais}. This contradicts that $\pi_1(C^0{\rm Aut}(\calL_p))$ is finite \cite[Lemma~3.1]{HJ}. Thus \smash{$\tp^{\rm W}_{\frac{3}{2}}$} is not identically zero.
 \end{proof}

\begin{Remark}
 We sketch a quicker proof of the main Theorem~\ref{thm:main} when $\om$ is K\"ahler (and is not $3$-Sasakian or $\mathbb{CP}^n$). By \cite[Corollary~8.1.19]{BG}, in this case $G=\ism(M,g)$ coincides with the group of strict contactomorphisms of $M$, where we only consider the identity components. Therefore, elements of $G $ commute with the circle of isometries given by the flow of the Reeb vector field (denoted by $a^I$ in Section~\ref{section3}). Thus
 $a^I$ lies in the center $Z$ of $G$.

From the fibration $Z\to G \to G/Z$ and using $\pi_2(G/Z) = 1$ (since $G$ is compact and connected), $\pi_1(Z)$ injects into $\pi_1(G)$. Since
$\bigl[a^I\bigr]$ is easily an element of infinite order in $\pi_1(Z)$, $\bigl[a^I\bigr]$ also has infinite order in $\pi_1(\ism(M,g))$.
\end{Remark}

\appendix

\section{The proof of Proposition~\ref{prop:temp}}\label{app:oldB}

\subsection{Pullbacks of forms}

 For $f\colon M\to N$ a differentiable map between finite-dimensional manifolds and $\omega\in \Omega^s(N)$, we have $ d_{M}f^*\omega = f^* d_N \omega$.
In \cite[Section~33.15]{KM}, this is
 extended to infinite dimensions in a very general setting which includes Fr\'echet manifolds.
 We give an alternative proof involving less notation for smooth Banach manifolds like $LM$.

 On an infinite-dimensional smooth manifold $N$, the exterior derivative can only be defined by the Cartan formula
\begin{align*} {\rm d}_N\omega\bigl(x^0,\dots X_{s}\bigr)_p={}& \sum_i (-1)^i X_i \bigl( \omega\bigl(x^0,\dots,\widehat{X_i},\dots, X_{s}\bigr)\bigr)\\
&{}{+}\, \sum_{i<j} (-1)^{i+j}\omega\bigl([X_i,X_j], x^0,\dots, \widehat{X_i},\dots, \widehat{X_j},\dots,X_{s}\bigr),
\end{align*}
where $X_i\in T_pN$ are extended to vector fields near $p$ using a chart map (see, e.g., \cite[Section~33.12]{KM}).
\begin{Lemma} \label{lem:2.1} Let $f\colon M\to N$ be a smooth map between smooth Banach manifolds, and let $\omega\in \Lambda^*(N)$.
Then $ {\rm d}_Mf^*\omega = f^* {\rm d}_N \omega$. In particular,
\[
{\rm d}_{\itbm}F^{L,*}\ocs^{\rm W}_{\kk} = F^{L,*} {\rm d}_{L\bm} \ocs^{\rm W}_{\kk}.
\]
\end{Lemma}

 \begin{proof} First assume that $f$ is an immersion on a neighborhood $U_p$ of a fixed $p\in M$.
For fixed vector fields $Y_i$ on $U_p$, set
 $g\colon f(U_p)\to \R$, $g(n) = \omega(f_*Y_1,\dots, f_*Y_s)_{n}$.
We have $(g\circ f)(m) = \omega(f_*Y_1,\dots, f_*Y_s)_{f(m)}$. Thus the identity $X_m(g\circ f) = (f_*X)
_{f(m)}(g)$
becomes
\[X_m(\omega(f_*Y_1,\dots, f_*Y_s)) = (f_*X)_{f(m)}(\omega(f_*Y_1,\dots, f_*Y_s))\]
for $\omega \in \Lambda^s(N)$.
Dropping $m, f(m)$, we get
\begin{align*}f^* {\rm d}_N\omega(X_0,\dots, X_{s})={}& {\rm d}_N\omega(f_*X_0,\dots, f_*X_{s})\\
={}&\sum_i (-1)^i f_*X_i \bigl( \omega\bigl(f_*X_0,\dots,\widehat{f_*X_i},\dots, f_*X_{s}\bigr)\bigr)\\
&
{+}\, \sum_{i<j} (-1)^{i+j}\omega\bigl([f_*X_i,f_*X_j], f_*X_0,\dots, \widehat{f_*X_i},\dots, \widehat{f_*X_j},\dots,f_*X_{s}\bigr)\\
={}& \sum_i (-1)^i X_i \bigl( \omega\bigl(f_*X_0,\dots,\widehat{f_*X_i},\dots, f_*X_{s}\bigr)\bigr)\\
&
{+}\, \sum_{i<j} (-1)^{i+j}\omega\bigl(f_*[X_i,X_j], f_*X_0,\dots, \widehat{f_*X_i},\dots, \widehat{f_*X_j},\dots,f_*X_{s}\bigr)\\
={}& {\rm d}_M \omega (f_*X_0,\dots, f_*X_s) = {\rm d}_Mf^*\omega(X_0,\dots,X_s),
\end{align*}
where we use $[f_*X_i,f_*X_j] = f_*[X_i,X_j]$ for immersions.

In general,
consider the graph $G\colon M\to M\times N$, $G(m) = (m, f(m))$. Then $\pi_N\circ G = f$ for the projection $\pi_N\colon M\times N\to N$. (We similarly define $\pi_M$.)
$G$ is an immersion, with $G_*(Y) = (Y, f_*Y)$ taking a vector field on $M$ to a well-defined vector field on $M\times N$.

Fix $(m_0, n_0)\in M\times N$, and set $i_M\colon M\to N\times N$, $i_M\colon N\to M\times N$ by
$i_M(m) = (m, n_0)$, $i_N(n) = (m_0,n)$. If a vector $(X_0,Y_0)\in T_{(m_0,n_0)}M\times N$
is extended to a nearby vector field $(X,Y)$ with $X$ constant in $N$ directions and $Y$ constant in $M$ directions,
it is straightforward to apply the Cartan formula to derive the standard equality (usually abbreviated
$ {\rm d}_{M\times N} = {\rm d}_M + {\rm d}_N$)
\[ {\rm d}_{M\times N}\alpha_{(m_0,n_0)} = \pi_M^*[ {\rm d}_M (i_M^*\alpha)_{m_0}] + \pi_N^*[ {\rm d}_N (i_N^*\alpha)_{n_0}] \]
for $\alpha\in \Lambda^*(M\times N)$.
Since $\pi i_M\colon m\mapsto n_0$ (so ${\rm d}_Mi_M^*\pi^*\omega = 0$) and $\pi i_N = {\rm id}$, the argument above for the immersion $G$ yields
\begin{align*} {\rm d}_Mf^*\omega &= {\rm d}_M G^*\pi_N^*\omega = G^* {\rm d}_{M\times N}\pi_N^*\omega
= G^*[\pi_M^* {\rm d}_M i_M^*\pi^*\omega + \pi_N^* {\rm d}_N i_N^*\pi^*\omega] \\
&= G^*\pi_N^* {\rm d}_N i_N^*\pi^*\omega = f^* {\rm d}_N \omega.\tag*{\qed}
\end{align*}
\renewcommand{\qed}{}
\end{proof}

As usual, ${\rm d}_Mf^* = f^*{\rm d}_N$ gives the induced map $f^*\colon H^*(N,\R) \to H^*(M,\R)$ on de Rham cohomology.

\subsection{Local coordinates expression}

We work in local coordinates
$\bigl(x^0, x\bigr) = \bigl(x^0, x^1,\dots, x^{2k-1}\bigr)$ on $\itbm$.
Recall from (\ref{K_tensor}) that%
\begin{align*}
\MoveEqLeft{ K_{\nu \lambda_1 \dots \lambda_{2k-1}} }
= \sum_{\sigma } \sgn({\sigma})
R_{\lambda_{\sigma (1)} e_1 \nu}{}^{e_{2}}
R_{\lambda_{\sigma (2) } \lambda_{\sigma (3)} e_{3}}{}^{e_1}
R_{\lambda_{\sigma (4) } \lambda_{\sigma (5)} e_{4}}{}^{e_{3}}
\cdots
R_{\lambda_{\sigma (2k-2) } \lambda_{\sigma (2k-1)} e_{2}}{}^{e_{k-1}},
\nonumber
\end{align*}
for $\sigma$ a permutation of $\{1,\dots,2k-1\}$, and where $R_{ijk}^{\ \ \ \ell}$ are the components of the curvature tensor of the metric on a general manifold $M$. Then,
\begin{equation}\label{K}(2k-1)! K_{\nu \lambda_1 \dots \lambda_{2k-1}}{\rm d}x^\nu\otimes {\rm d}x^{\lambda_1}\wedge\cdots\wedge {\rm d}x^{\lambda_{2k-1}}
\end{equation}
is the local expression of an element of
$\Omega^1( {\bm}) \otimes \Omega^{2k-2} (\bm)$.
For $\gamma\in L\bm$ and $X_{\gamma,i}\in T_\gamma L\bm$, we have by (\ref{csg})
\begin{equation}
\label{CSW}
\ocs^{\rm W} (\gamma ) (X_{\gamma, 1} ,\dots, X_{\gamma , 2k-1}) \\
= \frac{k}{2^{k-2}}\int_0^{2\pi}
 K_{\nu \lambda_1 \dots \lambda_{2k-1}} (\gamma (\theta))
{\dot \gamma}^{\nu} (\theta )
X_{\gamma ,1}^{\lambda_{1}}(\theta ) \cdots X_{\gamma , 2k-1}^{\lambda_{2k-1}}
 {\rm d}\theta.
\end{equation}
Then \smash{$\ocs^{\rm W}
\in \Omega^{2k-1} (L{\bm})$}, since we have contracted out the $\nu$ index. Since the integrand in (\ref{CSW}) is tensorial, we can integrate over $[0, 2\pi]$ even if the image of $\gamma$ does not lie in one coordinate~chart.

From now on, we drop the constant $k/2^{k-2}$ in (\ref{CSW}).

\subsection[Computing F\^{}\{L,*\} d\_\{LM\} CS\^{}W]{Computing $\boldsymbol{F^{L,*}{\rm d}_{L\bm} \ocs^{\rm W}}$}

We have
\begin{align*}
 \MoveEqLeft{\bigl({\rm d}_{L\bm} \ocs^{\rm W}_\gamma \bigr)
 (X_{\gamma ,0}, X_{\gamma ,1},\dots, X_{\gamma ,2k-1} ) }\nonumber\\
 &=
 \sum_{a=0}^{2k-1}
 (-1)^a X_{\gamma ,a}
 \bigl(\ocs^{\rm W} \bigl(X_{\gamma ,0},\dots,
 \widehat{X_{\gamma ,a}},
 \dots, X_{\gamma ,2k-1}\bigr)\bigr) \\
&\qquad +
\sum_{a<b} (-1)^{a+b}
 \bigl(\ocs^{\rm W} \bigl( [X_{\gamma , a}, X_{\gamma , b}],
 X_{\gamma ,0},\dots, \widehat{X_{\gamma ,a}},
 \dots, \widehat{X_{\gamma , b}}, X_{\gamma ,2k-1}\bigr)\bigr) \nonumber \\
 &:= \sum_a (1)_a + \sum_{a<b}(2)_{a,b}\nonumber.
 \end{align*}
Let
$\gamma_s (\theta) \in L{\bm} $ be a family of loops with
$\gamma_0 (\theta ) = \gamma (\theta)$,
$({\rm d}/{\rm d}s)|_{s=0} \gamma_s = X_{\gamma , a} $.
Then
\begin{align}\label{eq:4}
& X_{\gamma ,a}
 \bigl(\ocs^{\rm W} \bigl(X_{\gamma ,0},\dots, \widehat{X_{\gamma ,a}},
 \dots, X_{\gamma ,2k-1}\bigr)\bigr) \nonumber \\
 &=
\int_{0}^{2\pi}
 {\frac{{\rm d}}{{\rm d}s}}\biggl|_{s=0}
\bigl[K_{\nu \lambda_0 \dots \widehat{\lambda_a}\dots \lambda_{2k-1}}
(\gamma_s (\theta )) {\dot{\gamma_s}^\nu}
X_{\gamma_s,0}^{\lambda_0}
\cdots \widehat{X_{\gamma_s ,a}^{\lambda_a}}\cdots
X_{\gamma_s,2k-1}^{\lambda_{2k-1}}
 {\rm d}\theta \bigr] \nonumber\\
&=
\int_{0}^{2\pi}
 \partial_{x^\mu} K_{\nu \lambda_0 \dots \widehat{\lambda_a}\dots\lambda_{2k-1}}
(\gamma(\theta))
X_{\gamma, a}^{\mu} \dot{\gamma}^\nu (\theta)
X_{\gamma ,0}^{\lambda_0 } (\theta)
\cdots \widehat{X_{\gamma ,a}^{\lambda_a}}(\theta)\cdots
X_{\gamma ,2k-1}^{\lambda_{2k-1} } (\theta)
 {\rm d}\theta \nonumber\\
&\quad{} +
\int_{0}^{2\pi}
K_{\nu \lambda_0 \dots \widehat{\lambda_a}\dots \lambda_{2k-1}} (\gamma (\theta) )
\dot{X}_{\gamma ,a}^\nu (\theta )
X_{\gamma ,0}^{\lambda_0 } (\theta)
\cdots\widehat{X_{\gamma ,a}^{\lambda_a}}(\theta)\cdots
X_{\gamma ,2k-1}^{\lambda_{2k-1} } (\theta) \\
&\quad{} +
\int_{0}^{2\pi}
K_{\nu \lambda_0 \dots\widehat{\lambda_a}\dots \lambda_{2k-1}} (\gamma (\theta) )
\dot{\gamma}^\nu (\theta)
\bigl(\delta_{X_{\gamma, a}}X_{\gamma, 0}^{\lambda_0}\bigr)(\theta)
X_{\gamma ,1}^{\lambda_1 } (\theta)
\cdots\widehat{X_{\gamma ,a}^{\lambda_a}}(\theta) \cdots X_{\gamma ,2k-1}^{\lambda_{2k-1} } (\theta) + \cdots\nonumber \\
&\quad{} +
\int_{0}^{2\pi}
K_{\nu \lambda_0 \dots\widehat{\lambda_a}\dots \lambda_{2k-1}} (\gamma (\theta) )
\dot{\gamma}^\nu (\theta)
X_{\gamma ,0}^{\lambda_0 } (\theta)
\cdots\widehat{X_{\gamma ,a}^{\lambda_a}}(\theta)\cdots
X_{\gamma,2k-2}^{\lambda_{2k-2} } (\theta) \nonumber
 \cdot
 \bigl(\delta_{X_{\gamma, a}}X_{\gamma, 2k-1}^{\lambda_{2k-1}} \bigr)(\theta).\nonumber
\end{align}
Here $\delta_{X_{\gamma,a}}$ is the Fr\'echet derivative or directional derivative in the direction of $X_{\gamma,a}$, so, e.g., in finite dimensions we have
\begin{equation}\label{eq:lie}[X,Y]^i = \delta_XY^i - \delta_YX^i.
\end{equation}
 Since the flow of a vector field $X_\gamma(\theta)$ on $LM$ is computed by the flow at each fixed $\theta$, the Lie derivative of two vector fields on $LM$ is computed at each $\theta$.
Thus (\ref{eq:lie}) is valid at a fixed $\theta$. Denote the last three lines of (\ref{eq:4}) by (\ref{eq:4})${}_a$. It follows that
\begin{equation*}
 \sum_{a=0}^{2k-1} (-1)^a ({\rm A.4})_a
 + \sum_{a<b} (2)_{a,b} =0
\end{equation*}
Therefore,
\begin{align*}
&\bigl( {\rm d}_{L\bm} \ocs^{\rm W}_\gamma \bigr)
 (X_{\gamma ,0}, X_{\gamma ,1},\dots, X_{\gamma ,2k-1} ) \\
&\quad= \sum_{a=0}^{2k-1} (-1)^a
\int_{0}^{2\pi}
 \partial_{x^\mu} K_{\nu \lambda_0 \dots \widehat{\lambda_a}\dots\lambda_{2k-1}}
(\gamma(\theta))
X_{\gamma, a}^{\mu} \dot{\gamma}^\nu (\theta)
X_{\gamma ,0}^{\lambda_0 } (\theta)
\cdots \widehat{X_{\gamma ,a}^{\lambda_a}}(\theta)\cdots
X_{\gamma ,2k-1}^{\lambda_{2k-1} } (\theta)
 {\rm d}\theta \nonumber\\
&\qquad{} +\sum_{a=0}^{2k-1} (-1)^a
\int_{0}^{2\pi}
K_{\nu \lambda_0 \dots \widehat{\lambda_a}\dots \lambda_{2k-1}} (\gamma (\theta) )
\dot{X}_{\gamma ,a}^\nu (\theta )
X_{\gamma ,0}^{\lambda_0 } (\theta)
\cdots\widehat{X_{\gamma ,a}^{\lambda_a}}(\theta) \cdots
X_{\gamma ,2k-1}^{\lambda_{2k-1} } (\theta) {\rm d}\theta.
\end{align*}

For the pullback, we consider
\smash{$\bigl( F^{L,*} d_{L{\bm}} \ocs^{\rm W}\bigr)
\bigl(\partial_{x^0}, \partial_{x^1},
\dots, \partial_{x^{2k-1}}\bigr)$}
as a function on $[0,1]\times U$, where \smash{$(U,x)=\bigl(x^1,\dots, x^{2k-1}\bigr)$} is a coordinate chart on $\bm$.
Then
\begin{align}\label{eq:pullback}
&\bigl( F^{L,*} {\rm d}_{L{\bm}} \ocs^{\rm W}\bigr)
(\partial_{x^0}, \partial_{x^1},
\dots, \partial_{x^{2k-1}})_{(x^0,x)} \\
&=
 {\rm d}_{L{\bm}} \ocs^{\rm W}
 \bigl( F^L_* \partial_{x^0}, F^L_* \partial_{x^1},
\dots, F^L_* \partial_{x^{2k-1}} \bigr)_{F(x^0,x)} \nonumber\\
&=
 {\rm d}_{L{\bm} } \ocs^{\rm W} \biggl(
\frac{\partial F^{\lambda_0}}{\partial x^0} \partial_{x^{\lambda_0}},
\frac{\partial F^{\lambda_1}}{\partial x^1}\partial_{x^{\lambda_1}},
\dots ,
\frac{\partial F^{\lambda_{2k-1}}}{\partial x^{2k-1}}\partial_{x^{\lambda_{2k-1}}}\biggr)_{F(x^0,x)}
\nonumber\\
&=
\sum_{a=0}^{2k-1} (-1)^a
\int_0^{2\pi} \partial_{x^\mu}
K_{\nu \lambda_0 \dots\widehat{\lambda_a}\dots \lambda_{2k-1}}
\bigl(F\bigl(x^0,\theta, x\bigr)\bigr)
\frac{\partial F^\mu}{\partial x^a}
\frac{\partial F^\nu}{\partial \theta}
\frac{\partial F^{\lambda_0}}{\partial x^0}\cdots
\widehat{\frac{\partial F^{\lambda_a}}{\partial x^a} }
\cdots \frac{\partial F^{\lambda_{2k-1}}}{\partial x^{2k-1}}
 {\rm d}\theta \nonumber\\
&\quad{}+
\sum_{a=0}^{2k-1} (-1)^a
\int_0^{2\pi}
K_{\nu \lambda_0 \dots \widehat{\lambda_a} \dots \lambda_{2k-1}}
\bigl(F\bigl(x^0, \theta , x\bigr)\bigr)
\frac{\partial^2 F^\nu}{\partial x^a \partial \theta}
\frac{\partial F^{\lambda_0}}{\partial x^0}
\cdots
\widehat{\frac{\partial F^{\lambda_a}}{\partial x^a}}
\cdots \frac{\partial F^{\lambda_{2k-1}}}{\partial x^{2k-1}}
 {\rm d}\theta.\nonumber
\end{align}
One term in the last equation in (\ref{eq:pullback}) vanishes. The proof is in Appendix~\ref{app:B}.

\begin{Lemma} \label{lem:2.2a}
\begin{align*}
\int_0^{2\pi} \sum_{a=0}^{2k-1} (-1)^a
\partial_{x^\mu}
K_{\nu \lambda_0 \dots \widehat{\lambda_a}\dots\lambda_{2k-1}}
\bigl(F\bigl(x^0,\theta, x\bigr)\bigr)
\frac{\partial F^\nu}{\partial \theta}
\frac{\partial F^\mu}{\partial x^a}
\frac{\partial F^{\lambda_0}}{\partial x^0}\cdots
\widehat{\frac{\partial F^{\lambda_a}}{\partial x^a} }
\cdots
\frac{\partial F^{\lambda_{2k-1}}}{\partial x^{2k-1}}
 {\rm d}\theta =0.
\end{align*}
\end{Lemma}

Thus, we have the following.

\begin{Proposition}
\begin{align}
\MoveEqLeft{F^{L,*} {\rm d}_{L{\bm}} \ocs^{\rm W}
(\partial_{x^0}, \partial_{x^1},
\dots, \partial_{x^{2k-1}})_{(x^0,x)}}\nonumber\\
&=
\sum_{a=0}^{2k-1} (-1)^a
\int_0^{2\pi}
K_{\nu \lambda_0 \dots \widehat{\lambda_a} \dots \lambda_{2k-1}}
\bigl(F\bigl(x^0, \theta , x\bigr)\bigr)
\frac{\partial^2 F^\nu}{\partial x^a \partial \theta}
\frac{\partial F^{\lambda_0}}{\partial x^0}
\widehat{\frac{\partial F^{\lambda_a}}{\partial x^a}}
\cdots
\frac{\partial F^{\lambda_{2k-1}}}{\partial x^{2k-1}}
 {\rm d}\theta.\!\!\!\label{eq:four}
\end{align}
\end{Proposition}

\subsection{Homotopies of loops of diffeomorphisms}

We now make the assumption that
\begin{equation}\label{eq:diff}
F\bigl(x^0, \theta,\cdot\bigr)\colon\ \bm\to\bm \ \text{is a diffeomorphism for all} \ \bigl(x^0,\theta\bigr)\in [0,1]\times S^1.
\end{equation}
Then \smash{$\bigl\{F_*\bigl(\partial/\partial x^i\bigr)\bigr\}_{i=1}^{2k-1}$} is a basis of {\smash{$T_{F(x^0,\theta,x)}\bm$}} for all $\bigl(x^0,\theta,x\bigr)$.
Therefore, there exist functions $\alpha^i =\alpha^i\bigl(x^0,\theta,x\bigr)$, $i = 1,\dots,2k-1$, such that
\begin{equation}\label{eq:five} F_*\biggl(\frac{\partial}{\partial x^0}\biggr) =
\alpha^iF_*\biggl(\frac{\partial}{\partial x^i}\biggr).
\end{equation}
Using coordinates $y^i = y^i\bigl(x^0,\theta, x\bigr)$ near $y = F\bigl(x^0,\theta,x\bigr)$, we have
\[F_*\biggl(\frac{\partial}{\partial x^0}\biggl|_{_{(x^0, \theta,x)}}\biggr) = \frac{\partial F^\lambda}{\partial x^0} \frac{\partial}{\partial y^\lambda}\biggl|_{_{y}} \in T_y\bm,\qquad
F_*\biggl(\frac{\partial}{\partial x^i}\biggl|_{_{(x^0, \theta,x)}}\biggr) = \frac{\partial F^\lambda}{\partial x^i} \frac{\partial}{\partial y^\lambda}\biggl|_{_{y}} \in T_y\bm.\]
Thus
\begin{equation}\label{eq:five1}
\frac{\partial F^\lambda}{\partial x^0} = \alpha^i\frac{\partial F^\lambda}{\partial x^i},\qquad
\frac{\partial^2 F^\lambda}{\partial\theta \partial x^0} =
\frac{\partial \alpha^i}{\partial\theta} \frac{\partial F^\lambda}{\partial x^i}
+ \alpha^i\frac{\partial^2 F^\lambda}{\partial\theta \partial x^i}.
\end{equation}
Plugging (\ref{eq:five1}) into (\ref{eq:four}) gives
\begin{align*}
\MoveEqLeft{F^{L,*} {\rm d}_{L{\bm}} CS^{\rm W}
(\partial_{x^0}, \partial_{x^1},
\dots, \partial_{x^{2k-1}})}\nonumber\\
&=
 \int_0^{2\pi} K_{\nu \lambda_1 \dots \lambda_{2k-1}}
\biggl( \frac{\partial \alpha^i}{\partial\theta} \frac{\partial F^\nu}{\partial x^i}
+ \alpha^i\frac{\partial^2 F^\nu}{\partial\theta \partial x^i}\biggr)
\frac{\partial F^{\lambda_1}}{\partial x^1}
\cdots
\frac{\partial F^{\lambda_{2k-1}}}{\partial x^{2k-1}}
 {\rm d}\theta\\
&\quad +\sum_{a=1}^{2k-1}
(-1)^a
\int_0^{2\pi}
K_{\nu \lambda_0 \dots \widehat{\lambda_a} \dots \lambda_{2k-1}}
\frac{\partial^2 F^\nu}{\partial x^a \partial \theta} \biggl( \alpha^i\frac{\partial F^{\lambda_0}}{\partial x^i}
\biggr) \frac{\partial F^{\lambda_1}}{\partial x^1}
\cdots
\widehat{\frac{\partial F^{\lambda_a}}{\partial x^a}}
\cdots
\frac{\partial F^{\lambda_{2k-1}}}{\partial x^{2k-1}}
 {\rm d}\theta. \nonumber
\end{align*}
The sum of the terms with the second partial derivatives vanishes.
\begin{Lemma}\label{lem:2.3a}
\begin{align*}
 0={}& \int_0^{2\pi} K_{\nu \lambda_1\dots\lambda_{2k-1}}
 \alpha^i\frac{\partial^2 F^\nu}{\partial\theta \partial x^i}
\frac{\partial F^{\lambda_1}}{\partial x^1}
\cdots
\frac{\partial F^{\lambda_{2k-1}}}{\partial x^{2k-1}}
 {\rm d}\theta\nonumber \\
&{+}\, \sum_{a=1}^{2k-1}
(-1)^a
\int_0^{2\pi}
K_{\nu \lambda_0 \dots \widehat{\lambda_a} \dots \lambda_{2k-1}}
\frac{\partial^2 F^\nu}{\partial x^a \partial \theta} \biggl( \alpha^i\frac{\partial F^{\lambda_0}}{\partial x^i}
\biggr) \frac{\partial F^{\lambda_1}}{\partial x^1}
\cdots
\widehat{\frac{\partial F^{\lambda_a}}{\partial x^a}}
\cdots
\frac{\partial F^{\lambda_{2k-1}}}{\partial x^{2k-1}}
 {\rm d}\theta.\nonumber
 \end{align*}
 \end{Lemma}

 The proof is in Appendix~\ref{app:B}.
 Changing the index $\nu$ to $\lambda_0$, we have proved the following.
 \begin{Lemma}\label{lem:2.2}
 Under assumption \eqref{eq:diff}, we have
 \begin{equation*}
 F^{L,*} {\rm d}_{L{\bm}} \ocs^{\rm W} (\partial_{x^0}, \partial_{x^1},\dots, \partial_{x^{2k-1}}) =
 \int_0^{2\pi} K_{\lambda_0 \lambda_1\dots\lambda_{2k-1}}
\frac{\partial \alpha^i}{\partial\theta} \frac{\partial F^{\lambda_0}}{\partial x^i}
\frac{\partial F^{\lambda_1}}{\partial x^1}
\cdots
\frac{\partial F^{\lambda_{2k-1}}}{\partial x^{2k-1}}
 {\rm d}\theta.
 \end{equation*}
 \end{Lemma}

 \subsection{Homotopies by loops of isometries}
We now make the further assumption that
\begin{equation}\label{eq:isom}
 F^I\bigl(x^0,\theta\bigr) := F\bigl(x^0, \theta,\cdot\bigr)\colon\ \bm\to\bm \ \text{is an isometry for all} \ \bigl(x^0,\theta\bigr)\in [0,1]\times S^1.
\end{equation}

The following computation finishes the proof of Proposition~\ref{customprop:app}.
\begin{Lemma}\label{lem:2.3}
 Under the assumption \eqref{eq:isom}, we have
 \[F^{L,*} {\rm d}_{L{\bm}} \ocs^{\rm W}=0.\]
\end{Lemma}

\begin{proof}
 By Lemma~\ref{lem:2.2}, at a fixed $x^0$, we have
 \begin{align*}
 \MoveEqLeft{ F^{L,*} {\rm d}_{L{\bm}} \ocs^{\rm W} (\partial_{x^0}, \partial_{x^1},\dots, \partial_{x^{2k-1}})|_{x} }\\
 &=
 \int_0^{2\pi} K_{\lambda_0 \lambda_1\dots\lambda_{2k-1}}\bigl(F\bigl(x^0,\theta,x\bigr)\bigr)
\frac{\partial \alpha^i}{\partial\theta} \frac{\partial F^{I,\lambda_0}}{\partial x^i}
\frac{\partial F^{I,\lambda_1}}{\partial x^1}
\cdots
\frac{\partial F^{I,\lambda_{2k-1}}}{\partial x^{2k-1}}
 {\rm d}\theta\\
 &=\int_0^{2\pi}\frac{\partial \alpha^i}{\partial\theta} K_{i1\dots 2k-1}(x) {\rm d}\theta
 = K_{i1\dots 2k-1}(x)\int_0^{2\pi}\frac{\partial \alpha^i}{\partial\theta} {\rm d}\theta=0,
 \end{align*}
 using (\ref{eq:11}).
 As in (\ref{CSW}), the integration over $[0, 2\pi]$ is valid, because the $\alpha^i$ are the components of a tensor/vector (\ref{eq:five}).
\end{proof}

The crucial identity \smash{$K_{i_0 i_1\dots i_{2k-1}}(x) =\bigl(F^{I,*}K\bigr){}_{i_0 i_1\dots i_{2k-1}}(x)$}, which holds for isometries and~fails for diffeomorphisms in general, is generalized in~\cite{MRTVI} to study the fundamental group of other finite and infinite-dimensional transformation groups.

\section{Proofs of Lemmas~\ref{lem:2.2a} and~\ref{lem:2.3a}}\label{app:B}

\begin{proof}[Proof of Lemma~\ref{lem:2.2a}]
Set $\dim (M) = 2k-1. $ Fix
 $x\in \bm$ and $\xi\in T_x\bm$.
 For $X_0, X_1,\dots$, $ X_{2k-1}\in T_x\bm$, set
 \begin{align}\label{eq:tk}\tilde K(X_0,\dots, X_{2k-1})_x
&=
\sum_{a=0}^{2k-1} (-1)^a\partial_{\lambda^\mu} K_{\nu \lambda_0 \dots \widehat{\lambda_a}\dots\lambda_{2k-1}}
(x)\xi^\nu
X_a^{\mu} X_0^{\lambda_0}\cdots \widehat{X_\mu^{\lambda_a}} \cdots X_{2k-1}^{\lambda_{2k-1}}.
\end{align}
If we show that the right-hand side of (\ref{eq:tk}) is skew-symmetric in $X_0,\dots, X_{2k-1}$, then
$\tilde K
$ is a~$2k$-form on $\bm$ and hence must vanish. The lemma follows by replacing $x$ with
$F\bigl(x^0,\theta,x\bigr)$,~$\xi$~with \smash{$({\rm d}/{\rm d}\theta)F\bigl(x^0,\theta, x\bigr)$}, and \smash{$X_i^{\lambda_i}$} with $\partial F^{\lambda_i}/\partial x^i$.

To check skew-symmetry in $X_0, X_1$, we write
\begin{gather}
{\tilde{K}(X_0,X_1, X_2,\dots, X_{2k-1})
}\nonumber\\
\quad{}= (\partial_{\lambda_0} K_{\nu \lambda_1\dots \lambda_{2k-1}}
- \partial_{\lambda_1} K_{\nu \lambda_0 \lambda_2 \dots\lambda_{2k-1}})
\xi^\nu
X_0^{\lambda_0} X_1^{\lambda_1} \cdots X^{2k-1}_{\lambda_{2k-1}}
\label{16}\\
\qquad{}+ (\partial_{\lambda_2} K_{\nu \lambda_0 \lambda_1 \widehat{\lambda_2} \lambda_3\dots \lambda_{2k-1}}
-\partial_{\lambda_3} K_{\nu \lambda_0 \lambda_1 \lambda_2\widehat{\lambda_3}\lambda_4\dots \lambda_{2k-1}}
+\dots - \partial_{\lambda_{2k-1}}K_{\nu\lambda_0\dots\lambda_{2k-1}})\nonumber\\
\qquad\quad{} \times \xi^\nu X_0^{\lambda_0} X_1^{\lambda_1} X_2^{\lambda_2}\cdots X_{2k-1}^{\lambda_{2k-1}},\label{17}\\
{\tilde K(X_1,X_0, X_2, \dots,X_{2k-1}) }\nonumber\\
\quad{}= (\partial_{\lambda_1} K_{\nu \lambda_0\lambda_2\lambda_3\dots\lambda_{2k-1}}
- \partial_{\lambda_0} K_{\nu \lambda_1 \lambda_2\dots\lambda_{2k-1}}) \xi^\nu
X_1^{\lambda_1} X_0^{\lambda_0} X_2^{\lambda_2}X_3^{\lambda_3}\cdots X_{2k-1}^{\lambda_{2k-1}} \label{18}\\
\qquad{}+ (\partial_{\lambda_2} K_{\nu \lambda_1 \lambda_0 \widehat{\lambda_2} \lambda_3\dots \lambda_{2k-1}}
-\partial_{\lambda_3} K_{\nu \lambda_1 \lambda_0 \lambda_2\widehat{\lambda_3}\lambda_4\dots \lambda_{2k-1}}
+\cdots - \partial_{\lambda_{2k-1}}K_{\nu\lambda_1 \lambda_0 \lambda_2\dots\lambda_{2k-1}})\!\!\nonumber\\
\qquad\quad{} \times \xi^\nu X_1^{\lambda_1} X_0^{\lambda_0} X_2^{\lambda_2}X_3^{\lambda_3}\cdots X_{2k-1}^{\lambda_{2k-1}}.\label{19}
\end{gather}
Then $\text{(\ref{16})} = - \text{(\ref{18})}$ by inspection, and $\text{(\ref{17})} = - \text{(\ref{19})}$, because $K$ is skew-symmetric in $\lambda_1,\dots,\lambda_{2k-1}$ by (\ref{K}).

We now check skew-symmetry in $X_1$, $X_2$, with all other cases being similar. We have
\begin{align}
\MoveEqLeft{\tilde K(X_0,X_2, X_1, X_3, X_4, \dots,X_{2k-1}) }\nonumber\\
&= \partial_{\lambda_0} K_{\nu \lambda_2\lambda_1\lambda_3\lambda_4\dots\lambda_{2k-1}}
\xi^\nu X_1^{\lambda_1} X_0^{\lambda_0} X_2^{\lambda_2}X_3^{\lambda_3}X_4^{\lambda_4} \cdots X_{2k-1}^{\lambda_{2k-1}}\label{20}\\
&\quad +\bigl( - \partial_{\lambda_2} K_{\nu \lambda_0\lambda_1\lambda_3\lambda_4\dots\lambda_{2k-1}}
+ \partial_{\lambda_1} K_{\nu \lambda_0\lambda_2\lambda_3\lambda_4\dots\lambda_{2k-1}} \bigr)\nonumber\\
&\qquad \qquad \times
\xi^\nu X_1^{\lambda_1} X_0^{\lambda_0} X_2^{\lambda_2}X_3^{\lambda_3}X_4^{\lambda_4} \cdots X_{2k-1}^{\lambda_{2k-1}}\label{21}\\
&\qquad +\bigl(- \partial_{\lambda_3} K_{\nu \lambda_0\lambda_2\lambda_1\lambda_4\dots\lambda_{2k-1}}
+ \partial_{\lambda_4} K_{\nu \lambda_0\lambda_2\lambda_1\lambda_3\lambda_5\dots\lambda_{2k-1}}
+\cdots\nonumber\\
&\qquad\qquad - \partial_{\lambda_{2k-1}} K_{\nu \lambda_0\lambda_2\lambda_1\lambda_3\lambda_4
\dots \lambda_{2k-2}}\bigr)\nonumber\\
&\qquad \qquad \times\xi^\nu X_0^{\lambda_0} X_2^{\lambda_2}X_1^{\lambda_1}X_3^{\lambda_3}X_4^{\lambda_4}\cdots X_{2k-1}^{\lambda_{2k-1}}.\label{22}
\end{align}
Then $\tilde K(X_0,X_2, X_1, X_3, \dots, X_{2k-1}) = - \tilde K(X_0,X_1, X_2, X_3, \dots, X_{2k-1}) $, because
(i) the skew-symmetry of $K$ implies the skew-symmetry of (\ref{20}) and the skew-symmetry of the three lines of (\ref{22}) in $\lambda_1$, $\lambda_2$; (ii) the two lines of
(\ref{21}) are explicitly skew-symmetric in $\lambda_1$, $\lambda_2$.
\end{proof}

\begin{proof}[Proof of Lemma~\ref{lem:2.3a}]
The terms with second partial derivatives are
\begin{align}
&K_{\nu\lambda_1\lambda_2\lambda_3\lambda_4\dots\lambda_{2k-1}}\biggl(\alpha^1 \frac{\partial^2F^\nu}{\partial x^1\partial\theta}+ \alpha^2 \frac{\partial^2F^\nu}{\partial x^2\partial\theta}
+\alpha^3 \frac{\partial^2F^\nu}{\partial x^3\partial\theta} +\cdots
+\alpha^{2k-1} \frac{\partial^2F^\nu}{\partial x^{2k-1}\partial\theta}
\biggr)\nonumber\\
&\qquad \quad{} \times \frac{\partial F^{\lambda_1}}{\partial x^1}\frac{\partial F^{\lambda_2}}{\partial x^2}
\frac{\partial F^{\lambda_3}}{\partial x^3}
\cdots\frac{\partial F^{\lambda_{2k-1}}}{\partial x^{2k-1}}\label{23}\\
&\qquad{} -K_{\nu\lambda_0\lambda_2\lambda_3\lambda_4\dots \lambda_{2k-1}}\frac{\partial^2F^\nu}{\partial x^1\partial\theta}
\biggl( \alpha^i\frac{\partial F^{\lambda_0}}{\partial x^i}\biggr)
\frac{\partial F^{\lambda_2}}{\partial x^2}
\frac{\partial F^{\lambda_3}}{\partial x^3}\frac{\partial F^{\lambda_4}}{\partial x^4}\cdots
\frac{\partial F^{\lambda_{2k-1}}}{\partial x^{2k-1}}
\label{24}\\
&\quad + K_{\nu\lambda_0\lambda_1\lambda_3\lambda_4\dots\lambda_{2k-1}}\frac{\partial^2F^\nu}{\partial x^2\partial\theta}
\biggl( \alpha^i\frac{\partial F^{\lambda_0}}{\partial x^i}\biggr)
\frac{\partial F^{\lambda_1}}{\partial x^1}\frac{\partial F^{\lambda_3}}{\partial x^3}
\frac{\partial F^{\lambda_4}}{\partial x^4}\cdots \frac{\partial F^{\lambda_{2k-1}}}{\partial x^{2k-1}}\label{25}\\
&\quad -\cdots\nonumber\\
&\quad - K_{\nu\lambda_0\lambda_1\lambda_2\lambda_3\lambda_4\dots\lambda_{2k-2}}\frac{\partial^2F^\nu}{\partial x^{2k-1}\partial\theta}
\biggl( \alpha^i\frac{\partial F^{\lambda_0}}{\partial x^i}\biggr)
\frac{\partial F^{\lambda_1}}{\partial x^1}\frac{\partial F^{\lambda_2}}{\partial x^2}
\frac{\partial F^{\lambda_3}}{\partial x^3}\frac{\partial F^{\lambda_4}}{\partial x^4}\cdots
\frac{\partial F^{\lambda_{2k-2}}}{\partial x^{2k-2}}.
\label{28}
\end{align}
In (\ref{24}), in the term $\alpha^i \bigl(\partial F^{\lambda_0}/\partial x^i\bigr)$, only the term $\alpha^1 \bigl(\partial F^{\lambda_0}/\partial x^1\bigr)$ is nonzero: for example, the~term
\[K_{\nu\lambda_0\lambda_2\lambda_3\lambda_4\dots \lambda_{2k-1}}\frac{\partial^2F^\nu}{\partial x^1\partial\theta}
\biggl( \alpha^2\frac{\partial F^{\lambda_0}}{\partial x^2}\biggr) \frac{\partial F^{\lambda_2}}{\partial x^2}
\frac{\partial F^{\lambda_3}}{\partial x^3}\frac{\partial F^{\lambda_4}}{\partial x^4}\cdots
\frac{\partial F^{\lambda_{2k-1}}}{\partial x^{2k-1}}\]
is skew-symmetric in $\lambda_0$, $\lambda_2$, and so vanishes. For the same reasons, the terms with
$\partial F^{\lambda_0}/\partial x^3$, $ \partial F^{\lambda_0}/\partial x^4, \dots, \partial F^{\lambda_0}/\partial x^{2k-1}$ vanish.
 Similarly, in (\ref{25}) only $\alpha^2 \bigl(\partial F^{\lambda_0}/\partial x^2\bigr)$ is nonzero, etc.,
and in (\ref{28}) only $\alpha^{2k-1} \bigl(\partial F^{\lambda_0}/\partial x^{2k-1}\bigr)$ is nonzero.

Thus (\ref{23})--(\ref{28}) becomes
\begin{align}
&K_{\nu\lambda_1\lambda_2\lambda_3\lambda_4\dots\lambda_{2k-1}}\biggl(\alpha^1 \frac{\partial^2F^\nu}{\partial x^1\partial\theta}+ \alpha^2 \frac{\partial^2F^\nu}{\partial x^2\partial\theta}
+\alpha^3 \frac{\partial^2F^\nu}{\partial x^3\partial\theta} +\cdots
+\alpha^{2k-1} \frac{\partial^2F^\nu}{\partial x^{2k-1}\partial\theta}
\biggr)\nonumber\\
&\quad\qquad{} \times \frac{\partial F^{\lambda_1}}{\partial x^1}\frac{\partial F^{\lambda_2}}{\partial x^2}
\frac{\partial F^{\lambda_3}}{\partial x^3}
\cdots\frac{\partial F^{\lambda_{2k-1}}}{\partial x^{2k-1}}\label{29}\\
&\qquad{}-K_{\nu\lambda_0\lambda_2\lambda_3\lambda_4\dots \lambda_{2k-1}}\frac{\partial^2F^\nu}{\partial x^1\partial\theta}
\biggl( \alpha^1\frac{\partial F^{\lambda_0}}{\partial x^1}\biggr)
\frac{\partial F^{\lambda_2}}{\partial x^2}
\frac{\partial F^{\lambda_3}}{\partial x^3}\frac{\partial F^{\lambda_4}}{\partial x^4}\cdots
\frac{\partial F^{\lambda_{2k-1}}}{\partial x^{2k-1}}
\label{30}\\
&\qquad{} + K_{\nu\lambda_0\lambda_1\lambda_3\lambda_4\dots\lambda_{2k-1}}\frac{\partial^2F^\nu}{\partial x^2\partial\theta}
\biggl( \alpha^2\frac{\partial F^{\lambda_0}}{\partial x^2}\biggr)
\frac{\partial F^{\lambda_1}}{\partial x^1}\frac{\partial F^{\lambda_3}}{\partial x^3}
\frac{\partial F^{\lambda_4}}{\partial x^4}\cdots \frac{\partial F^{\lambda_{2k-1}}}{\partial x^{2k-1}}
\label{31}\\
&\qquad -\cdots\nonumber\\
&\qquad{} - K_{\nu\lambda_0\lambda_1\lambda_2\lambda_3\lambda_4\dots\lambda_{2k-2}}\frac{\partial^2F^\nu}{\partial x^{2k-1}\partial\theta}
\biggl( \alpha^{2k-1}\frac{\partial F^{\lambda_0}}{\partial x^{2k-1}}\biggr)\nonumber\\
&\quad\qquad{}\times
\frac{\partial F^{\lambda_1}}{\partial x^1}\frac{\partial F^{\lambda_2}}{\partial x^2}
\frac{\partial F^{\lambda_3}}{\partial x^3}\frac{\partial F^{\lambda_4}}{\partial x^4}\cdots
\frac{\partial F^{\lambda_{2k-2}}}{\partial x^{2k-2}}.
\label{34}
\end{align}

If we replace $\lambda_0$ in (\ref{30}) with $\lambda_1$, then the term
\[K_{\nu\lambda_1\lambda_2\lambda_3\lambda_4\dots\lambda_{2k-1}}\biggl(\alpha^1 \frac{\partial^2F^\nu}{\partial x^1\partial\theta}\biggr)\frac{\partial F^{\lambda_1}}{\partial x^1}\frac{\partial F^{\lambda_2}}{\partial x^2}
\frac{\partial F^{\lambda_3}}{\partial x^3}\frac{\partial F^{\lambda_4}}{\partial x^4}\cdots
\frac{\partial F^{\lambda_{2k-1}}}{\partial x^{2k-1}}\]
in (\ref{29}) cancels with (\ref{30}).
If we replace $\lambda_0$ in (\ref{31}) with $\lambda_2$, then the term
\[K_{\nu\lambda_1\lambda_2\lambda_3\lambda_4\dots\lambda_{2k-1}}\biggl(\alpha^2 \frac{\partial^2F^\nu}{\partial x^2\partial\theta}\biggr)\frac{\partial F^{\lambda_1}}{\partial x^1}\frac{\partial F^{\lambda_2}}{\partial x^2}
\frac{\partial F^{\lambda_3}}{\partial x^3}\frac{\partial F^{\lambda_4}}{\partial x^4}\cdots
\frac{\partial F^{\lambda_{2k-1}}}{\partial x^{2k-1}}\]
in (\ref{29}) cancels with (\ref{31}).
Continuing, if we replace $\lambda_0$ in (\ref{34}) with $\lambda_{2k-1}$, then the term
\[K_{\nu\lambda_1\lambda_2\lambda_3\lambda_4\dots\lambda_{2k-1}}\biggl(\alpha^{2k-1} \frac{\partial^2F^\nu}{\partial x^{2k-1}\partial\theta}\biggr)\frac{\partial F^{\lambda_1}}{\partial x^1}\frac{\partial F^{\lambda_2}}{\partial x^2}
\frac{\partial F^{\lambda_3}}{\partial x^3}\frac{\partial F^{\lambda_4}}{\partial x^4}\cdots
\frac{\partial F^{\lambda_{2k-1}}}{\partial x^{2k-1}}\]
in (\ref{29}) cancels with (\ref{34}).

Thus (\ref{29})--(\ref{34}) sum to zero, which proves the lemma.
\end{proof}

\section{Proof of Claim~\ref{claim1}}\label{app:claim}

\noindent {\bf Claim~\ref{claim1}:} For each $m\in M$,
\begin{equation}\label{Seqn1}
S_{4n+1,4n+2}(m) =
(-1)^{n+1} 2^{2n+1}(2n+1)
\sum_{\sigma \in \mathfrak S_{4n}} \sgn(\sigma) J_{\sigma_1 \sigma_2} \cdots J_{\sigma_{4n-1} \sigma_{4n}}(m),
\end{equation}
{\it where $\mathfrak S_{4n}$ is the permutation group of $\{1,\dots, 4n\}$.}

\begin{proof}
By (\ref{eq:cyan4}) and (\ref{eq:11aa}), $S_{4n+1,4n+2}$ is the term in $K_{1,12\dots 4n+1}$ with power $p^{4n+2}$. In the expression for $K_{1,12\dots 4n+1}$ in (\ref{K_tensor}) with $\nu=1$, we plug in Lemma~\ref{lem1}\,(i)--(iv) for the curvature terms.
Since $\nu=1$, the only way to obtain a term with power $p^{4n+2}$
is if a permutation $\sigma\in \mathfrak{C}_{4n+1}$ on the right-hand side of (\ref{K_tensor}) has $\sigma_1 :=\sigma(1) =1$. See Lemma~\ref{lem1}\,(iii) in particular.

However, to get the indices to agree with (\ref{Seqn1}), we need an adjustment. While in (\ref{CSWlocal}) we took local coordinates $\bigl(x^1,x^2,\dots, x^{4n+1}\bigr)$ on $M$ with $\partial_{x^1} = \xi$, we now take coordinates $\bigl(x^0, x^1,\dots, x^{4n}\bigr)$ such that
$\partial_{x^0} = \xi$. In particular, $K_{1, 1 2\dots 4n+1}$ in (\ref{eq:cyan4}) is replaced by
$K_{0, 01 2\dots 4n}$. The only permutations in $\mathfrak{C}_{4n+1}$, which is now the permutation group of $\{0,1,2,\dots,4n\}$, which give power $p^{4n+2}$ have $\sigma_0=0$. The set of such permutations equals $\mathfrak{C}_{4n}$, the permutations of~$\{1,2,\dots,4n\}$. Therefore, by Lemma~\ref{lem1}\,(iii), we get
\begin{align*}
S_{4n+1,4n+2}p^{4n+2}
={}& \sum_{\sigma\in \mathfrak{C}_{4n}}
 \sgn(\sigma)
\barcurvature{c_1 a_1 c_2}{a_{n}}
\xi^{c_1} \xi^{c_2}
\barcurvature{\sigma_1 \sigma_2 b_1}{a_1}
\barcurvature{\sigma_3 \sigma_4 a_2}{b_1}
\barcurvature{\sigma_5 \sigma_6 b_2}{a_2}
\barcurvature{\sigma_7 \sigma_8 a_3}{b_2}
\cdots \\
 &\times
\barcurvature{\sigma_{4n-3} \sigma_{4n-2} b_n}{a_{n-1}}
\barcurvature{\sigma_{4n-1} \sigma_{4n} a_{n}}{b_n}
\pmod{p^{4n+1}}
\\
\equiv{}& \sum_{\sigma\in \mathfrak{C}_{4n}}
 \sgn(\sigma)
\bigl(-p^2 \delta_{a_1}^{a_{n}} \bigr)
\bigl(\barcurvature{\sigma_1 \sigma_2 b_1}{a_1}
\barcurvature{\sigma_3 \sigma_4 a_2}{b_1}\bigr)
\bigl(\barcurvature{\sigma_5 \sigma_6 b_2}{a_2}
\barcurvature{\sigma_7 \sigma_8 a_3}{b_2}\bigr)
\cdots \\
 &\times
\bigl(\barcurvature{\sigma_{4n-3} \sigma_{4n-2} b_{n}}{a_{n-1}}
\barcurvature{\sigma_{4n-1} \sigma_{4n} a_n}{b_{n}} \bigr)
\pmod{p^{4n+1}},
 \end{align*}
where
$\bigl({\rm mod}\ p^{4n+1}\bigr)$ mods out all polynomials of degree at most $4n+1$, i.e.,
keeps only terms with power $p^{4n+2}$. (We have relabeled $e_1,\dots,e_{2n}$ in
(\ref{K_tensor}) for $2k-1 = 4n+1$ with $a_1,\dots,a_n, b_1,\dots, b_n$.)

By (\ref{curvconv}) and Lemma~\ref{lem1}\,(i), this becomes
\begin{gather*}
S_{4n+1,4n+2}p^{4n+2}
 \equiv -p^{4n+2} \sum_{\sigma\in \mathfrak{C}_{4n}}
 \sgn(\sigma)
 \delta_{a_1}^{a_n}
 \bigl({A_1^\prime}\bigr)_{\sigma_1 \sigma_2 \sigma_3 \sigma_4 a_2}{}^{a_1}
 \bigl({A_2^\prime}\bigr)_{\sigma_5 \sigma_6 \sigma_7 \sigma_8 a_3}{}^{a_2}\cdots
 \\
\hphantom{S_{4n+1,4n+2}p^{4n+2} \equiv -p^{4n+2} \sum_{\sigma\in \mathfrak{C}_{4n}}}{}
 {}\times
\bigl({A_n^\prime}\bigr)_{\sigma_{4n-3} \sigma_{4n-2} \sigma_{4n-1} \sigma_{4n} a_n}{}^{a_{n-1}}
\pmod{p^{4n+1}},
\end{gather*}
where
\begin{gather*}
{\bigl({A_1^\prime}\bigr)_{\sigma_1 \sigma_2 \sigma_3 \sigma_4 a_2}{}^{a_1} = \barcurvature{\sigma_1 \sigma_2 b_1}{a_1}
\barcurvature{\sigma_3 \sigma_4 a_2}{b_1} }
 \\ \nonumber
\quad{}=\bigl[-J_{\sigma_2 b_1} J_{\sigma_1}{}^{a_1}
 +J_{\sigma_1 b_1} J_{\sigma_2}{}^{a_1}
 +2 J_{\sigma_1 \sigma_2} J_{b_1}{}^{a_1} \bigr]
\bigl[-J_{\sigma_4 a_2} J_{\sigma_3}{}^{b_1} +J_{\sigma_3 a_2} J_{\sigma_4}{}^{b_1}
 +2 J_{\sigma_3 \sigma_4} J_{a_2}{}^{b_1}\bigr],
\\ \nonumber
{\bigl({A_2^\prime}\bigr)_{\sigma_5 \sigma_6 \sigma_7 \sigma_8 a_3}{}{}^{a_2}
= \barcurvature{\sigma_5 \sigma_6 b_2}{a_2}
\barcurvature{\sigma_7 \sigma_8 a_3}{b_2}}
 \\ \nonumber
\quad{}=\bigl[-J_{\sigma_6 b_2} J_{\sigma_5}{}^{a_2}
 +J_{\sigma_5 b_2} J_{\sigma_6}{}^{a_2}
 +2 J_{\sigma_5 \sigma_6} J_{b_2}{}^{a_2} \bigr]
 \bigl[-J_{\sigma_8 a_3} J_{\sigma_7}{}^{b_2}
 +J_{\sigma_7 a_3} J_{\sigma_8}{}^{b_2}
 +2 J_{\sigma_7 \sigma_8} J_{a_3}{}^{b_2} \bigr],
 \\ \nonumber
\quad \ \, {} \vdots \\
{\bigl({A_n^\prime}\bigr)_{\sigma_{4n-3} \sigma_{4n-2} \sigma_{4n-1} \sigma_{4n} a_n}{}^{a_{n-1}}
=\barcurvature{\sigma_{4n-3} \sigma_{4n-2} b_{n}}{a_{n-1}}
\barcurvature{\sigma_{4n-1} \sigma_{4n} a_n}{b_{n}}}
\\ \nonumber
\quad{}=\bigl[-J_{\sigma_{4n-2} b_{n}} J_{\sigma_{4n-3}}{}^{a_{n-1}}
 +J_{\sigma_{4n-3} b_{n}} J_{\sigma_{4n-2}}{}^{a_{n-1}}
 +2 J_{\sigma_{4n-3} \sigma_{4n-2}} J_{b_{n}}{}^{a_{n-1}} \bigr]
 \\ \nonumber
 \quad\quad{} \times
\bigl[-J_{\sigma_{4n} a_n} J_{\sigma_{4n-1}}{}^{b_{n}}
 +J_{\sigma_{4n-1} a_n}
 J_{\sigma_{4n}}{}^{b_{n}}
 +2 J_{\sigma_{4n-1} \sigma_{4n}}
 J_{a_n}{}^{b_n} \bigr].
\end{gather*}
We expand out $\bigl(A_1^\prime\bigr)$,
\begin{gather*}
{\bigl({A_1^\prime}\bigr)_{\sigma_1 \sigma_2 \sigma_3 \sigma_4 a_2}{}^{a_1} }
= J_{\sigma_2 b_1} J_{\sigma_1}{}^{a_1} J_{\sigma_4 a_2} J_{\sigma_3}{}^{b_1}
 -J_{\sigma_2 b_1} J_{\sigma_1}{}^{a_1} J_{\sigma_3 a_2} J_{\sigma_4}{}^{b_1}
 -2 J_{\sigma_2 b_1} J_{\sigma_1}{}^{a_1} J_{\sigma_3 \sigma_4} J_{a_2}{}^{b_1} \\ \nonumber
\qquad{}
 -J_{\sigma_1 b_1} J_{\sigma_2}{}^{a_1} J_{\sigma_4 a_2} J_{\sigma_3}{}^{b_1}
 +J_{\sigma_1 b_1} J_{\sigma_2}{}^{a_1} J_{\sigma_3 a_2} J_{\sigma_4}{}^{b_1}
 +2 J_{\sigma_1 b_1} J_{\sigma_2}{}^{a_1} J_{\sigma_3 \sigma_4} J_{a_2}{}^{b_1}
\\ \nonumber
\qquad{}
 -2J_{\sigma_1 \sigma_2} J_{b_1}{}^{a_1} J_{\sigma_4 a_2} J_{\sigma_3}{}^{b_1}
 +2J_{\sigma_1 \sigma_2} J_{b_1}{}^{a_1} J_{\sigma_3 a_2} J_{\sigma_4}{}^{b_1}
 +4 J_{\sigma_1 \sigma_2} J_{b_1}{}^{a_1} J_{\sigma_3 \sigma_4} J_{a_2}{}^{b_1 }.
 \end{gather*}
Since, e.g.,
\begin{equation}\label{eq:minus}
J_{\sigma_2 b_1}J_{\sigma_3}{}^{b_1} = J^{a}_{b_1} g_{a\sigma_2}J_{\sigma_3}{}^{b_1}
= -\delta^a_{\sigma_3} g_{a\sigma_2} =- g_{\sigma_2 \sigma_3},
\end{equation}
$\bigl(A_1^\prime\bigr)$ reduces to
\begin{align}
\bigl({A_1^\prime}\bigr)_{\sigma_1 \sigma_2 \sigma_3 \sigma_4 a_2}{}^{a_1}
={}& {-}J_{\sigma_1}^{\ a_1} J_{\sigma_4 a_2} g_{\sigma_2 \sigma_3}
+J_{\sigma_1}{}^{a_1} J_{\sigma_3 a_2} g_{\sigma_2 \sigma_4}
+2 J_{\sigma_1}{}^{a_1} J_{\sigma_3 \sigma_4} g_{\sigma_2 a_2}
 \nonumber\\
&{+}\, J_{\sigma_2}{}^{a_1} J_{\sigma_4 a_2} g_{\sigma_1 \sigma_3}
- J_{\sigma_2}{}^{a_1} J_{\sigma_3 a_2} g_{\sigma_1 \sigma_4}
-2 J_{\sigma_2}{}^{a_1} J_{\sigma_3 \sigma_4} g_{\sigma_1 a_2}\nonumber
\\
&{+}\, 2J_{\sigma_1 \sigma_2} J_{\sigma_4 a_2}\delta_{\sigma_3}{}^{a_1}
-2J_{\sigma_1 \sigma_2} J_{\sigma_3 a_2}\delta_{\sigma_4}{}^{a_1}
-4J_{\sigma_1 \sigma_2} J_{\sigma_3 \sigma_4}\delta_{a_2}{}^{a_1}.\label{Y3}
\end{align}
Terms in (\ref{Y3}) containing $g_{\sigma_i\sigma_j}$ (as opposed to terms containing
$g_{\sigma_i a_j}$) do not contribute to $S_{4n+1,4n+2}$. Indeed, by the symmetry of $g$,
for fixed $\sigma$ the term in $S_{4n+1,4n+2}$ formally of the form
$\sgn(\sigma) J\cdot J\cdots J\cdot g_{\sigma_i\sigma_j}$ is
 cancelled by the term with
$ (ij)\sigma$ in cycle notation.

As a result, we have
\begin{align}
S_{4n+1,4n+2}p^{4n+2}
 = p^{4n+2} \sum_{\sigma_0 = 0}& \sgn(\sigma) \delta_{a_1}^{a_n}
 ({A_1})_{\sigma_1 \sigma_2 \sigma_3 \sigma_4 a_2}{}^{a_1}
 \bigl({A_2^\prime}\bigr)_{\sigma_5 \sigma_6 \sigma_7 \sigma_8 a_3}{}^{a_2}\cdots\nonumber \\
&\times
\bigl({A_n^\prime}\bigr)_{\sigma_{4n-3} \sigma_{4n-2} \sigma_{4n-1} \sigma_{4n} a_n}{}^{a_{n-1}},\label{3.5add5}
\end{align}
with
\begin{align*}
({A_1})_{\sigma_1 \sigma_2 \sigma_3 \sigma_4 a_2}{}^{a_1}
={}&2J_{\sigma_1}{}^{a_1} J_{\sigma_3 \sigma_4} g_{\sigma_2 a_2}
-2 J_{\sigma_2}{}^{a_1} J_{\sigma_3 \sigma_4} g_{\sigma_1 a_2}
+2J_{\sigma_1 \sigma_2} J_{\sigma_4 a_2} \delta_{\sigma_3}{}^{a_1}
\\ \nonumber
&
{-}\,2J_{\sigma_1 \sigma_2} J_{\sigma_3 a_2} \delta_{\sigma_4}{}^{a_1}
-4J_{\sigma_1 \sigma_2} J_{\sigma_3 \sigma_4} \delta_{a_2}{}^{a_1}
\\ \nonumber
={}&{-}4 J_{\sigma_3 \sigma_4} J_{\sigma_1}{}^{a_1} g_{\sigma_2 a_2}
 -4 J_{\sigma_1 \sigma_2} J_{\sigma_3 a_2} \delta_{\sigma_4}{}^{a_1}
 -4J_{\sigma_1 \sigma_2} J_{\sigma_3 \sigma_4} \delta_{a_2}{}^{a_1}
\\ \nonumber
={}& {-}2^2 \bigl[J_{\sigma_1 \sigma_2} J_{\sigma_3 a_2} \delta_{\sigma_4}{}^{a_1}
 + J_{\sigma_3 \sigma_4} J_{\sigma_1}{}^{a_1} g_{\sigma_2 a_2}
 + J_{\sigma_1 \sigma_2} J_{\sigma_3 \sigma_4} \delta_{a_2}{}^{a_1}\bigr]
\\ \nonumber
={}&{-}2^2 J_{\sigma_1 \sigma_2}
 \bigl(J_{\sigma_3 a_2} \delta_{\sigma_4}{}^{a_1}
 +J_{\sigma_3}{}^{a_1} g_{\sigma_4 a_2}
 + J_{\sigma_3 \sigma_4} \delta_{a_2}{}^{a_1} \bigr).
\end{align*}
(To obtain the third line, we replace
$2J_{\sigma_1}{}^{a_1} J_{\sigma_3 \sigma_4} g_{\sigma_2 a_2} $ by
$-2 J_{\sigma_2}{}^{a_1} J_{\sigma_3 \sigma_4} g_{\sigma_1 a_2} $, and
$2J_{\sigma_1 \sigma_2} J_{\sigma_4 a_2} \delta_{\sigma_3}{}^{a_1} $
by $-2J_{\sigma_1 \sigma_2} J_{\sigma_3 a_2} \delta_{\sigma_4}{}^{a_1}$, using the sign reversing ``change of variables" $\sigma\mapsto (12)\sigma$.
In the last line, we replaced $ J_{\sigma_3 \sigma_4} J_{\sigma_1}{}^{a_1} g_{\sigma_2 a_2}$ with
$J_{\sigma_1 \sigma_2} J_{\sigma_3}{}^{a_1} g_{\sigma_4 a_2} $ using the sign preserving change of variables
$\sigma\mapsto (13)(24)\sigma$. Strictly speaking, these substitutions are valid only after we plug $(A_1)$ back into~(\ref{3.5add5}).)

Doing the same computations for $\bigl(A_2^\prime\bigr),\dots, \bigl(A_n^\prime \bigr)$, we get
\begin{align}\label{3.5add7}
&{S_{4n+1,4n+2}} \\
 &\quad{}= \sum_{\sigma_0 = 0} \sgn(\sigma) (\delta_{a_1}^{a_n})
 ({A_1})_{\sigma_1 \sigma_2 \sigma_3 \sigma_4 a_2}{}^{a_1}
 ({A_2})_{\sigma_5 \sigma_6 \sigma_7 \sigma_8 a_3}{}^{a_2} \cdots
({A_n})_{\sigma_{4n-3} \sigma_{4n-2} \sigma_{4n-1} \sigma_{4n} a_n}{}^{a_{n-1}}, \nonumber
\end{align}
where
\begin{gather*} 
({A_2})_{\sigma_5 \sigma_6 \sigma_7 \sigma_8 a_3}{}^{a_2}
=
\bigl(-2^2\bigr) \bigl[J_{\sigma_5 \sigma_6}
 \bigl(J_{\sigma_7 a_3} \delta_{\sigma_8}{}^{a_2}
 +J_{\sigma_7}{}^{a_2} g_{\sigma_8 a_3}
 + J_{\sigma_7 \sigma_8} \delta_{a_3}{}^{a_2} \bigr) \bigr],
\\ \nonumber
\quad \ \, \vdots \\ \nonumber
 ({A_n})_{\sigma_{4n-3} \sigma_{4n-2} \sigma_{4n-1} \sigma_{4n} a_n} {}^{a_{n-1}} \nonumber\\
\qquad{}= \bigl(-2^2\bigr)
 \bigl[J_{\sigma_{4n-3} \sigma_{4n-2}} \bigl(J_{\sigma_{4n-1} a_n} \delta_{\sigma_{4n}}{}^{a_{n-1}}
 +J_{\sigma_{4n-1}}{}^{a_{n-1}} g_{\sigma_{4n} a_n}
+ J_{\sigma_{4n-1} \sigma_{4n}} \delta_{a_n}{}^{a_{n-1}} \bigr) \bigr]. \nonumber
\end{gather*}
We now begin to simplify (\ref{3.5add7}).
\begin{align*}
{S_{4n+1,4n+2}}
 = (-1) \sum_{\sigma_0 = 0}{}& \sgn(\sigma) (\delta_{a_1}^{a_n})
({A_{12}})_{\sigma_1 \sigma_2 \sigma_3 \sigma_4 \sigma_5 \sigma_6 \sigma_7 \sigma_8 a_3}{}^{a_1}
 ({A_3})_{\sigma_9 \sigma_{10} \sigma_{11} \sigma_{12} a_4}{}^{a_3}\cdots
\\ \nonumber
&
\times ({A_n})_{\sigma_{4n-3} \sigma_{4n-2} \sigma_{4n-1} \sigma_{4n} a_n}{}^{a_{n-1}},
\end{align*}
where
\begin{gather}
{\bigl({A_{12}}\bigr)_{\sigma_1 \sigma_2 \sigma_3 \sigma_4
\sigma_5 \sigma_6 \sigma_7 \sigma_8 a_3}{}^{a_1}
:= \bigl({A_1}\bigr)_{\sigma_1 \sigma_2 \sigma_3 \sigma_4 a_2}{}^{a_1}
 \bigl({A_2}\bigr)_{\sigma_5 \sigma_6 \sigma_7 \sigma_8 a_3}{}^{a_2} }\nonumber\\
\qquad{}=
\bigl(-2^2\bigr)^2 \bigl[
 J_{\sigma_1 \sigma_2 }
\bigl(J_{\sigma_3 a_2} \delta_{\sigma_4}{}^{a_1}
+ J_{\sigma_3 \sigma_4} J_{\sigma_3}^{a_1}g_{\sigma_4 a_1} +
J_{\sigma_3 \sigma_4} \delta_{a_2}{}^{a_1} \bigr)\bigr]
 \nonumber\\
 \qquad\quad{} \times
\bigl[J_{\sigma_5 \sigma_6} \bigl(J_{\sigma_7 a_3} \delta_{\sigma_8}{}^{a_2}
+ J_{\sigma_7}{}^{a_2} g_{\sigma_8 a_3}
+ J_{\sigma_7 \sigma_8} \delta_{a_2}{}^{a_3} \bigr)\bigr]
 \nonumber\\
\qquad{}=
\bigl(-2^2\bigr)^2\cdot J_{\sigma_1 \sigma_2} J_{\sigma_5 \sigma_6}
\bigl[J_{\sigma_3 a_2} \delta_{\sigma_4}{}^{a_1}
+J_{\sigma_3}{}^{a_1} g_{\sigma_4 a_2}
+J_{\sigma_3 \sigma_4} \delta_{a_2}{}^{a_1} \bigr]\nonumber\\
\quad\qquad{} \times
\bigl[J_{\sigma_7 a_3} \delta_{a_8}{}^{a_2}
+J_{\sigma_7}{}^{a_2} g_{\sigma_8 a_3}
+J_{\sigma_7 \sigma_8} \delta_{a_3}{}^{a_2} \bigr]
 \nonumber\\
\qquad{}=
\bigl(-2^2\bigr)^2\cdot
J_{\sigma_1 \sigma_2} J_{\sigma_5 \sigma_6}
\bigl[ J_{\sigma_3 a_2} \delta_{\sigma_4}{}^{a_1} J_{\sigma_7 a_3}
\delta_{\sigma_8}{}^{a_2}
+J_{\sigma_3 a_2}\delta_{\sigma_4}{}^{a_1} J_{\sigma_7}{}^{a_2}
g_{\sigma_8 a_3} \nonumber
\\
\qquad\qquad{} +J_{\sigma_3 a_2} \delta_{\sigma_4}{}^{a_1} J_{\sigma_7 \sigma_8}
\delta_{a_3}{}^{a_2}
+J_{\sigma_3}{}^{a_1} g_{\sigma_4 a_2} J_{\sigma_7 a_3}
\delta_{\sigma_8}{}^{a_2}
+J_{\sigma_3}{}^{a_1} g_{\sigma_4 a_2} J_{\sigma_7}{}^{a_2}
g_{\sigma_8 a_3} \nonumber \\
 \qquad\qquad{} +J_{\sigma_3}{}^{a_1} g_{\sigma_4 a_2} J_{\sigma_7 \sigma_8}
\delta_{a_3}{}^{a_2}
+J_{\sigma_3 \sigma_4} \delta_{a_2}{}^{a_1} J_{\sigma_7 a_3}
\delta_{\sigma_5}{}^{a_2}
+J_{\sigma_3 \sigma_4} \delta_{a_2}{}^{a_1} J_{\sigma_7}{}^{a_2}
g_{\sigma_8 a_3} \nonumber\\
\qquad\qquad{} +J_{\sigma_3 \sigma_4} \delta_{a_2}{}^{a_1} J_{\sigma_7 \sigma_8}
\delta_{a_3}{}^{a_2} \bigr]\nonumber\\
\qquad{} =
\bigl(-2^2\bigr)^2\cdot J_{\sigma_1 \sigma_2} J_{\sigma_5 \sigma_6}
\bigl[J_{\sigma_3 \sigma_8} J_{\sigma_7 a_3} \delta_{\sigma_4}{}^{a_1}
+J_{\sigma_3 a_3} \delta_{\sigma_4}{}^{a_1} J_{\sigma_7 \sigma_8}
J_{\sigma_3}{}^{a_1} J_{\sigma_7 \sigma_4} g_{\sigma_8 a_3}
 \nonumber\\
 \qquad\qquad{} +J_{\sigma_3}{}^{a_1} g_{\sigma_4 a_3} J_{\sigma_7 \sigma_8}
+J_{\sigma_3 \sigma_4} J_{\sigma_7 a_3} \delta_{\sigma_8}^{a_1}
+J_{\sigma_3 \sigma_4} J_{\sigma_7}{}^{a_1} g_{\sigma_8 a_3}
+J_{\sigma_3 \sigma_4} J_{\sigma_7 \sigma_8} \delta_{a_3}{}^{a_1}\bigr]
 \nonumber\\
\qquad{} =
\bigl(-2^2\bigr)^2\cdot J_{\sigma_1 \sigma_2} J_{\sigma_5 \sigma_6}
\bigl[-J_{\sigma_3 \sigma_4} J_{\sigma_7 a_3} \delta_{\sigma_8}{}^{a_1}
+J_{\sigma_7 a_3} \delta_{\sigma_8}{}^{a_1} J_{\sigma_3 \sigma_4}
-J_{\sigma_7}{}^{a_1} J_{\sigma_3 \sigma_4} g_{\sigma_8 a_3}
 \nonumber\\
 \qquad\qquad{} +J_{\sigma_7}{}^{a_1} g_{\sigma_8 a_3} J_{\sigma_3 \sigma_4}
+J_{\sigma_3 \sigma_4} J_{\sigma_7 a_3} \delta_{\sigma_8}^{a_1}
+J_{\sigma_3 \sigma_4} J_{\sigma_7}{}^{a_1} g_{\sigma_8 a_3}
+J_{\sigma_3 \sigma_4} J_{\sigma_7 \sigma_8} \delta_{a_3}{}^{a_1}\bigr]
 \nonumber\\
\qquad{}=
\bigl(-2^2\bigr)^2\cdot J_{\sigma_1 \sigma_2} J_{\sigma_5 \sigma_6}
\bigl[ J_{\sigma_3 \sigma_4}
\bigl(-J_{\sigma_7 a_3} \delta_{\sigma_8}{}^{a_1}
+J_{\sigma_7 a_3} \delta_{\sigma_8}{}^{a_1}
-J_{\sigma_7}{}^{a_1} g_{\sigma_8 a_3}
+J_{\sigma_7}{}^{a_1} g_{\sigma_8 a_3}
\nonumber\\
\qquad\qquad{} +J_{\sigma_7 a_3} \delta_{\sigma_8}{}^{a_1} +J_{\sigma_7}{}^{a_1} g_{\sigma_8 a_3}
+J_{\sigma_7 \sigma_8} \delta_{a_3}{}^{a_1}\bigr) \bigr]
\nonumber\\
\qquad{}= \bigl(-2^2\bigr)^2 \cdot J_{\sigma_1 \sigma_2} J_{\sigma_3 \sigma_4}
J_{\sigma_5 \sigma_6}
\bigl[J_{\sigma_7 a_3} \delta_{\sigma_8}{}^{a_1}
+J_{\sigma_7}{}^{a_1} g_{\sigma_8 a_3}
+J_{\sigma_7 \sigma_8} \delta_{a_3}{}^{a_1}\bigr].\label{3.5add10}
\end{gather}

Continuing to simplify (\ref{3.5add7}), we have
\begin{gather*}
 {(A_{123})_{\sigma_1\dots\sigma_{12}a_4}{}^{a_1}
:= (A_{12})_{\sigma_1 \dots \sigma_8 a_3}{}^{a_1}
 ({A_3})_{\sigma_9 \sigma_{10} \sigma_{11} \sigma_{12} a_4}{}^{a_3} }\\
\qquad{} = \bigl(-2^2\bigr)^2 \cdot J_{\sigma_1 \sigma_2} J_{\sigma_3 \sigma_4}J_{\sigma_5 \sigma_6}
\bigl[J_{\sigma_7 a_3} \delta_{\sigma_8}{}^{a_1}
+J_{\sigma_7}{}^{a_1} g_{\sigma_8 a_3}
+J_{\sigma_7 \sigma_8} \delta_{a_3}{}^{a_1}\bigr]\\
\qquad\quad{} \times \bigl(-2^2\bigr)\cdot J_{\sigma_9\sigma_{10}}\bigl[J_{\sigma_{11}a_4}\delta_{\sigma_{12}}{}^{a_3}
+J_{\sigma_{11}}{}^{a_3} g_{\sigma_{12} a_4} + J_{\sigma_{11}\sigma_{12}}\delta_{a_4}{}^{a_3}\bigr]\\
\qquad{}= \bigl(-2^2\bigr)^3\cdot J_{\sigma_1 \sigma_2} J_{\sigma_3 \sigma_4}J_{\sigma_5 \sigma_6} J_{\sigma_7 \sigma_8}
J_{\sigma_9 \sigma_{10}}
\bigl[J_{\sigma_{11}a_4}\delta_{\sigma_{12}}{}^{a_1}+ J_{\sigma_{11}}{}^{a_1}g_{\sigma_{12}a_4} + J_{\sigma_{11}\sigma_{12}}\delta_{a_4}{}^{a_1}\bigr],
\end{gather*}
where the last line follows from computations as in (\ref{3.5add10}).

In the end, we obtain
\begin{gather}
{S_{4n+1,4n+2}} =(-1)^{n+1} 2^{2n}
\sum_{\sigma_0=0} \sgn (\sigma)
J_{\sigma_1 \sigma_2} J_{\sigma_3 \sigma_4} \cdot\dots\cdot
J_{\sigma_{4n-3} \sigma_{4n-2}}
 \nonumber\\
 \hphantom{{S_{4n+1,4n+2}} =(-1)^{n+1} 2^{2n}
\sum_{\sigma_0=0}}{}
\times \bigl[J_{\sigma_{4n-1} a_1} \delta_{\sigma_{4n}}{}^{a_1}
+J_{\sigma_{4n-1}}{}^{a_1}g_{\sigma_{4n} a_{1}}+J_{\sigma_{4n-1} \sigma_{4n}} \delta_{a_1}{}^{a_1} \bigr]\nonumber\\
\hphantom{{S_{4n+1,4n+2}}}{}
=(-1)^{n+1}2^{2n+1}(2n+1)
\sum_{\sigma_0=0} \sgn (\sigma)
J_{\sigma_1 \sigma_2} J_{\sigma_3 \sigma_4} \cdots
J_{\sigma_{4n-1} \sigma_{4n}}.\label{3.5add12}
\end{gather}
This proves the claim.
\end{proof}

\section{Proof of Proposition~\ref{app1}}\label{app:E}

\noindent {\bf Proposition 5.2.} {\em Let $\pi\colon \omp\to \om$ be the fibration with $\om$ K\"ahler. For $\bxi$ the unit tangent vector to the fibers of }
$\pi$,
\begin{align*}\imath_{\bxi} a^{L,*}\ocs^{\rm W}_{4k+1,2} &=
(2k+1)2\cdot \pi^*\tr\bigl(\bom^{2k}\bigr)
=
(-1)^k (4k+2)(2\pi)^{2k+1}(2k)!
\cdot \pi^*\tilde p_k\bigl(\bom\bigr).
\end{align*}

\begin{proof}
Let $\bxi = e_1, e_2,\dots, e_{4n+1}$ be an orthonormal frame of $\omp$ at $ m$ with $\bxi$ tangent to the fiber of
the $S^1$ action and \smash{$\{e_i\}_{i=2}^{4n+1}$} a horizontal lift of an orthonormal frame $\{\eb_i = \pi_*e_i\}$ at $\bar m = \pi( m)$.
We must show that
\smash{$a^{L,*}\ocs{}^{\rm W}_{4k+1,2}(\bxi, e_2,\dots,e_{4k+1})$} at $\bar m$ is a specific multiple of \smash{$\tr\bigl(\bom^{2k}\bigr)(\pi_*e_2,\dots, \pi_*e_{4k+1})$} at $\pi(\bar m)$. (More precisely, we have to show this for every subset of $\{e_2,\dots, e_{4n+1}\}$ of size $4k$, but the case we treat carries over to all other cases.) In the proof, we will use Lemma~\ref{lem1} for curvature terms.

We denote $ e_r$ by $r$, so, e.g., $\br\bigl( e_{\sigma(2)}, e_{\sigma(3)}, e_{\ell_2}, e_{\ell_1}\bigr)$ is denoted $\br(\sigma(2), \sigma(3), \ell_2, \ell_1)$.
Since ${a^L_* e_r = e_r}$, as in
(\ref{csg}),
we have
\begin{align} \label{twob}
\MoveEqLeft{a^{L,*}\ocs^{\rm W}_{4k+1}(e_1,\dots,e_{4k+1}) }\nonumber\\
&=
\frac{2k+1}{2^{2k-1}}\sum_{\sigma\in \mathfrak{C}_{4k+1}} {\rm sgn}(\sigma) \int_{S^1}\tr\bigl[
 \bigl(R\bigl(e_{\sigma(1)},\cdot\bigr)\bxi\bigr)
 \bigl(\Omega^M\bigr)^{2k}\bigl(e_{\sigma(2)},\dots,e_{\sigma(4k+1)} \bigr)\bigr]\\
 &=
 \frac{2k+1}{2^{2k-1}}\sum_{\sigma\in \mathfrak{C}_{4k+1}} {\rm sgn}(\sigma) \int_{S^1}
 \br(\sigma(1), \ell_1,\bxi,r) \br(\sigma(2), \sigma(3), \ell_2, \ell_1)
 \br(\sigma(4), \sigma(5), \ell_3, \ell_2)\cdots\nonumber\\
 &\hphantom{=\frac{2k+1}{2^{2k-1}}\sum_{\sigma\in \mathfrak{C}_{4k+1}} }{}
 \times \br(\sigma(4k-2), \sigma(4k-1), \ell_{4k}, \ell_{4k-1})
 \br(\sigma(4k),\sigma(4k+1), r, \ell_{4k}).\nonumber
\end{align}

We want to compute the terms in (\ref{twob}) of order $p^2$.
These terms come from $(\alpha)$ permutations with $\sigma(1) = 1$, and ($\beta$) permutations with $\sigma(1)\neq 1$.

We claim the $(\beta)$ terms contribute zero \big(for all powers of $p^2$\big). The term
$ \br(\sigma(1), \ell_1,\bxi,r)$ with $\sigma(1)\neq 1$ is zero unless $\ell_1 = 1$ and $\sigma(1) = r$. Note that $\br(r, 1, \bxi, r) = -p^2$. Thus
\begin{gather*}
(\beta)=\frac{2k+1}{2^{2k-1}}\cdot p^2\ints \sum_{\sigma(1) = r\atop \ell_1 = 1} -\sgn(\sigma)
\br(\sigma(2), \sigma(3), \ell_2, 1 = \bxi)\cdots\\
\hphantom{(\beta)=\frac{2k+1}{2^{2k-1}}\cdot p^2\ints \sum_{\sigma(1) = r\atop \ell_1 = 1} }{} \times
 \br(\sigma(4k),\sigma(4k+1), r = \sigma(1), \ell_{4k}).\nonumber
 \end{gather*}
In the term $\br(\sigma(2), \sigma(3), \ell_2, 1 = \bxi)$, we
get zero unless either [$\sigma(2) = 1$ and $\ell_2 = \sigma(3)$] or [$\sigma(3) = 1$
and $\ell_2 = \sigma(2)]$. Therefore,
\begin{gather*}
(\beta)
=\frac{2k+1}{2^{2k-1}}\cdot p^2
\ints \sum_{\sigma(1) = r\atop \ell_1 = 1}\sum_{\sigma(2) = 1\atop \sigma(3) = \ell_2}\sgn(\sigma)
\br(\sigma(4), \sigma(5), \ell_3, \ell_2)\cdots \\
\hphantom{(\beta)=\frac{2k+1}{2^{2k-1}}\cdot p^2
\ints \sum_{\sigma(1) = r\atop \ell_1 = 1}\sum_{\sigma(2) = 1\atop \sigma(3) = \ell_2}}{}
\times
 \br(\sigma(4k),\sigma(4k+1), r = \sigma(1), \ell_{4k})\\
\hphantom{(\beta)=}{} - \frac{2k+1}{2^{2k-1}}\cdot p^2
\cdot \ints\sum_{\sigma(1) = r\atop \ell_1 = 1}\sum_{\sigma(3) = 1\atop \sigma(2) = \ell_2}\sgn(\sigma)
\br(\sigma(4), \sigma(5), \ell_3, \ell_2)\cdots\\
\hphantom{(\beta)=- \frac{2k+1}{2^{2k-1}}\cdot p^2\cdot \ints\sum_{\sigma(1) = r\atop \ell_1 = 1}\sum_{\sigma(3) = 1\atop \sigma(2) = \ell_2}}{}
\times
 \br(\sigma(4k),\sigma(4k+1), r = \sigma(1), \ell_{4k}).
 \end{gather*}

For fixed $\ell_2$, there is a bijection between $\{\sigma \mid \sigma(2) = 1,\, \sigma(3) = \ell_2\}$
and
$\{\tau\mid \tau(3) = 1,\, \tau(2) = \ell_2\}$ given by $\sigma\mapsto \tau =\sigma(1\ell_2)$ in cycle notation. Since $\sgn(\sigma) = -\sgn(\tau)$, we get
\begin{align}\nonumber
(\beta)
=\frac{2k+1}{2^{2k-1}}\cdot p^2
\ints \sum_{\sigma(1) = r\atop \ell_1 = 1}\sum_{\sigma(2) = 1\atop \sigma(3) = \ell_2}&\sgn(\sigma)
\br(\sigma(4), \sigma(5), \ell_3, \ell_2)\cdots\\
&\times
 \br(\sigma(4k),\sigma(4k+1), r = \sigma(1), \ell_{4k}).\label{beta2a}
 \end{align}
The last term in (\ref{beta2a}) is
\begin{equation}\label{temp}\br(\sigma(4k),\sigma(4k+1), r = \sigma(1), \ell_{4k}) =
\br(\sigma(4k),\sigma(4k+1), \bxi, \ell_{4k}).
\end{equation}
This term vanishes if $\ell_{4k} = \bxi$. If $\ell_{4k}\neq \bxi$, then since $\sigma(2) = 1$, we have $\sigma(4k)\neq \bxi$, ${\sigma(4k+1)\neq \bxi}$. Thus (\ref{temp}) vanishes in all cases.
Therefore, $(\beta) = 0$.

The $(\alpha)$ term is
\begin{align*}
\begin{split}
(\alpha)= \frac{2k+1}{2^{2k-1}}\cdot p^2\ints \sum_{\sigma(1) = 1\atop \ell_1 = r\neq 1}{}& \sgn(\sigma)
\br(\sigma(2), \sigma(3), \ell_2, \ell_1 = r)\cdots\\
& \times
 \br(\sigma(4k),\sigma(4k+1), \ell_1 =r, \ell_{4k}).\nonumber
\end{split}
 \end{align*}
 In Lemma~\ref{lem1},
 a nonzero product of terms of types (i) and (iii)
 with one term having $\bxi$ and having a power $p^2$ must
include exactly one term from (iii) and only the first term on the right-hand side of (i).
For $\bar r = \pi_* e_r = \bar e_r$, etc., we get
 \begin{align*}(\alpha)
&= \frac{(2k+1)2\pi}{2^{2k-1}}\cdot p^2
 \sum_{\sigma(1) = 1\atop \ell_1 = r\neq 1} \sgn(\sigma)
\overline{R}\bigl(\bar\sigma(2), \bar\sigma(3), \bar\ell_2, \bar\ell_1 = \bar r\bigr)\cdots \\
&\hphantom{= \frac{(2k+1)2\pi}{2^{2k-1}}\cdot p^2
 \sum_{\sigma(1) = 1\atop \ell_1 = r\neq 1}}{}
\times
 \overline{R}\bigl(\bar\sigma(4k),\bar\sigma(4k+1), \bar\ell_1 =\bar r, \bar \ell_{4k}\bigr)\nonumber\\
 &= \frac{2k+1}{2^{2k-1}}\cdot p^2 \cdot 2^{2k}2\pi\cdot \Tr\bigl(\bom^{2k}\bigr)(\bar e_2,\dots, \bar e_{4k+1}),\\
 &= (-1)^k (4k+2)(2\pi)^{2k}(2k)!\cdot p^2 \cdot \tilde p_k\bigl(\bom\bigr)
 (\bar e_2,\dots, \bar e_{4k+1})\\
 &= (-1)^k (4k+2)(2\pi)^{2k+1}(2k)!\cdot
p^2 \cdot \pi^* \tilde p_k\bigl(\bom\bigr)
 (e_2,\dots, e_{4k+1}).\tag*{\qed}
 \end{align*}\renewcommand{\qed}{}
\end{proof}

\section[Manifolds of dimension 4n+2]{Manifolds of dimension $\boldsymbol{4n+2}$}\label{appb}

We prove that $S_{4n+3, 4n+4}$ vanishes if dim $\om = 4n+2$. We proved the stronger result that \smash{$\ocs{}^{\rm W}_{3} = 0$} for any $3$-manifold \cite[Proposition~2.7]{MRT4}.

\begin{Proposition} For $\dim\om = 4n+2$, $S_{4n+3,4n+4}(\omp) = 0$.
\end{Proposition}

We do the case where dim $\om = 6$ to keep the notation down. The proof immediately extends to the general case.

\begin{proof}As in the proof of Claim~\ref{claim1} in Appendix~\ref{app:claim}, we have
\begin{align}\label{S78}
S_{7,8} &=\ssuzero\barcurvature{0\ell_1 0}{n}\barcurvature{\sigma_1 \sigma_2\ell_2}{\ell_1}
\barcurvature{\sigma_3\sigma_4\ell_3}{\ell_2}\barcurvature{\sigma_5\sigma_6 n}{\ell_3}\\
&= \ssuzero \bigl(p^2\delta_{a_1}^b\bigr)\cdot p^2\bigl[-J_{\sigma_1}^{a_1}J_{\sigma_2a_2} + J_{\sigma_1a_2}J_{\sigma_2}^{a_1}
+2J_{\sigma_1\sigma_2}J_{a_2}^{a_1}\bigr]\nonumber\\
&\hphantom{=\sum_{\sigma_0 = 0}}{}\times p^2\bigl[-J_{\sigma_3}^{a_2}J_{\sigma_4a_3} + J_{\sigma_3a_3}J_{\sigma_4}^{a_2}
+2J_{\sigma_3\sigma_4}J_{a_3}^{a_2}\bigr]
 \cdot p^2\bigl[-J_{\sigma_5}^{a_3}J_{\sigma_6b} + J_{\sigma_5 b}J_{\sigma_6}^{a_3}
+2J_{\sigma_5\sigma_6}J_{b}^{a_3}\bigr]\nonumber\\
&= p^8\ssuzero\bigl[-J_{\sigma_1}^{a_1}J_{\sigma_2a_2} + J_{\sigma_1a_2}J_{\sigma_2}^{a_1}
+2J_{\sigma_1\sigma_2}J_{a_2}^{a_1}\bigr]\nonumber\\
&\hphantom{= p^8\sum_{\sigma_0 = 0}} \times
\bigl[-J_{\sigma_3}^{a_2}J_{\sigma_4a_3} + J_{\sigma_3a_3}J_{\sigma_4}^{a_2}
+2J_{\sigma_3\sigma_4}J_{a_3}^{a_2}\bigr]
 \!\cdot \!\bigl[-J_{\sigma_5}^{a_3}J_{\sigma_6a_1} + J_{\sigma_5 a_1}J_{\sigma_6}^{a_3}
+2J_{\sigma_5\sigma_6}J_{a_1}^{a_3}\bigr]. \nonumber
\end{align}
Because $J_{ab}J^b_c
= -g_{ba}$ by (\ref{eq:minus}),
the product of the first
two expressions in square brackets simplifies to
\begin{gather*}
S_{7,8} = p^8\ssuzero\bigl[-2J_{\sigma_1}^{a_1}J_{\sigma_3\sigma_4}g_{\sigma_2a_3}
+2g_{\sigma_1a_3}J_{\sigma_2}^{a_1}J_{\sigma_3\sigma_4}
-2J_{\sigma_1\sigma_2}\bigl(-\delta_{\sigma_3}^{a_1}\bigr)J_{\sigma_4a_3}\\
\hphantom{S_{7,8} = p^8\ssuzero\bigl[}{}
+2J_{\sigma_3a_3}\bigl(-\delta_{\sigma_4}^{a_1}\bigr)J_{\sigma_1\sigma_2}
+4J_{\sigma_1\sigma_2}J_{\sigma_3\sigma_4}\bigl(-\delta_{a_3}^{a_1}\bigr)\bigr]\\
\hphantom{S_{7,8} = p^8\sum_{\sigma_0 = 0}\bigl[}{}
 \times \bigl[-J_{\sigma_5}^{a_3}J_{\sigma_6a_1} + J_{\sigma_5 a_1}J_{\sigma_6}^{a_3}
+2J_{\sigma_5\sigma_6}J_{a_1}^{a_3}\bigr].
\end{gather*}
Taking the product of the terms inside the first square brackets with the terms inside the second square brackets, we get 15 terms, all of which simplify. For example, the product of the first terms in each square brackets gives
\[2J_{\sigma_1}^{a_1}J_{\sigma_3\sigma_4}g_{\sigma_2a_3}J_{\sigma_5}^{a_3}J_{\sigma_6a_1}
= -2g_{\sigma_1\sigma_6}J_{\sigma_3\sigma_4}J_{\sigma_2\sigma_5}.\]
For a term with a Kronecker delta, we have
\begin{align*}-2J_{\sigma_1\sigma_2}\bigl(-\delta_{\sigma_3}^{a_1}\bigr)J_{\sigma_4a_3}\cdot J_{\sigma_5}^{a_3}J_{\sigma_6a_1}
&= -2J_{\sigma_1\sigma_2}J_{\sigma_4\sigma_3}J^{a_3}_{\sigma_5}J_{\sigma_6\sigma_3}
= -2J_{\sigma_1\sigma_2}J^k_{a_3}g_{k\sigma_4}J^{a_3}_{\sigma_5}J_{\sigma_6\sigma_3}\\
&= -2J_{\sigma_1\sigma_2}\bigl(-\delta^k_{\sigma_5}\bigr)g_{k\sigma_4}J_{\sigma_6\sigma_3}
= -2J_{\sigma_1\sigma_2}g_{\sigma_5\sigma_4}J_{\sigma_6\sigma_3}.
\end{align*}
Similarly, every product is of the form $g_{\sigma_i\sigma_j} J_{\sigma_k\sigma_\ell}J_{\sigma_m\sigma_n}$,
except for the product of the two last terms, which is
\[8J_{\sigma_1\sigma_2}J_{\sigma_3\sigma_4}\bigl(-\delta_{a_3}^{a_1}\bigr)J_{\sigma_5\sigma_6}J_{a_1}^{a_3}
=0,\]
since $\delta_{a_3}^{a_1}J_{a_1}^{a_3} = J_{a_1}^{a_1} = 0$.

In summary, every nonzero term in $S_{7,8}$ is of the form $\sgn(\sigma)g_{\sigma_i\sigma_j} J_{\sigma_k\sigma_\ell}J_{\sigma_m\sigma_n}$, where $\sigma = (i,j,k,\ell,m,n)\in \mathfrak S_6$. Under the change of variables $\sigma\mapsto \sigma(12)$, $\sigma$ changes sign, but the term
$g_{\sigma_i\sigma_j} J_{\sigma_k\sigma_\ell}J_{\sigma_m\sigma_n}$ does not change sign. Thus the terms corresponding to $\sigma$ and $\sigma(12)$ in $S_{7,8}$ cancel, so $S_{7,8} = 0$.
\end{proof}

\begin{Remark} In this proof, it was crucial that $S_{7,8}$ contains an odd number of bracket terms in (\ref{S78}); if there are an even number of terms, the cancellation in the previous paragraph does not occur. This is
where the hypothesis dim $\om = 4n+2$ is used.
\end{Remark}

\subsection*{Acknowledgements}

We would like to thank Yoshinobu Kamishima for helpful conversations and the referees for their insightful comments.

\pdfbookmark[1]{References}{ref}
\LastPageEnding


\begin{thebibliography}{99}
\footnotesize\itemsep=0pt

\bibitem{Blair}
Blair D.E., Riemannian geometry of contact and symplectic manifolds, 2nd ed., \textit{Progr. Math.}, Vol.~203,
 \href{https://doi.org/10.1007/978-0-8176-4959-3}{Birkh\"auser}, Boston, MA,
 2010.

\bibitem{BG}
Boyer C.P., Galicki K., Sasakian geometry, \textit{Oxford Math. Monogr.}, Oxford
 University Press, Oxford, 2008.

\bibitem{Bredon}
Bredon G.E., Introduction to compact transformation groups, \textit{Pure Appl.
 Math.}, Vol.~46, Academic Press, New York, 1972.

\bibitem{Egi1}
Egi S., Calculation of the coefficient of the highest power of~$p$
 ({T}heorem~3.4), 2020,
 \url{https://github.com/egisatoshi/EMR-Paper-Computation}.

\bibitem{Egi3}
Egi S., {T}hurston-calculation.pdf, 2020,
 \url{https://github.com/egison/egison/blob/master/sample/math/geometry/thurston.egi}.

\bibitem{fgls}
Fedosov B.V., Golse F., Leichtnam E., Schrohe E., The noncommutative residue
 for manifolds with boundary,
 \href{https://doi.org/10.1006/jfan.1996.0142}{\textit{J.~Funct. Anal.}}
 \textbf{142} (1996), 1--31.

\bibitem{HJ}
Hsiang W.C., Jahren B., A~note on the homotopy groups of the diffeomorphism
 groups of spherical space forms, in Algebraic {$K$}-Theory, {P}art~{II}
 ({O}berwolfach, 1980), \textit{Lecture Notes in Math.}, Vol.~967,
 \href{https://doi.org/10.1007/BFb0061901}{Springer}, Berlin, 1982, 132--145.

\bibitem{Klingenberg}
Klingenberg W., Lectures on closed geodesics, \textit{Grundlehren Math. Wiss.},
 Vol.~230, \href{https://doi.org/10.1007/978-3-642-61881-9}{Springer}, Berlin,
 1978.

\bibitem{Kodaira}
Kodaira K., On the structure of compact complex analytic surfaces.~{I},
 \href{https://doi.org/10.2307/2373157}{\textit{Amer.~J.~Math.}} \textbf{86}
 (1964), 751--798.

\bibitem{KM}
Kriegl A., Michor P.W., The convenient setting of global analysis,
 \textit{Math. Surveys Monogr.}, Vol.~53,
 \href{https://doi.org/10.1090/surv/053}{American Mathematical Society},
 Providence, RI, 1997.

\bibitem{alh-unpublished}
Larrain-Hubach A., Wodzicki--{C}hern classes and vortex equations,
 {u}npublished.

\bibitem{lrst}
Larrain-Hubach A., Rosenberg S., Scott S., Torres-Ardila F., Characteristic
 classes and zeroth order pseudodifferential operators, in Spectral Theory and
 Geometric Analysis, \textit{Contemp. Math.}, Vol.~535,
 \href{https://doi.org/10.1090/conm/535/10539}{American Mathematical Society},
 Providence, RI, 2011, 141--158,
 \href{http://arxiv.org/abs/1003.0067}{arXiv:1003.0067}.

\bibitem{MRTVI}
Maeda Y., Rosenberg S., The geometry of loop spaces~{V}: {F}undamental groups
 of geometric transformation groups,
 \href{http://arxiv.org/abs/2510.01566}{arXiv:2510.01566}.

\bibitem{MRT3}
Maeda Y., Rosenberg S., Torres-Ardila F., The geometry of loop spaces~{I}:
 {$H^s$}-{R}iemannian metrics,
 \href{https://doi.org/10.1142/S0129167X15400029}{\textit{Internat.~J.~Math.}}
 \textbf{26} (2015), 1540002, 26~pages,
 \href{http://arxiv.org/abs/1405.4231}{arXiv:1405.4231}.

\bibitem{MRT4}
Maeda Y., Rosenberg S., Torres-Ardila F., The geometry of loop spaces~{II}:
 {C}haracteristic classes,
 \href{https://doi.org/10.1016/j.aim.2015.10.001}{\textit{Adv. Math.}}
 \textbf{287} (2016), 485--518.

\bibitem{MRT4-cor}
Maeda Y., Rosenberg S., Torres-Ardila F., The geometry of loop spaces~{II}: {C}orrections,
 \href{http://arxiv.org/abs/2405.00651}{arXiv:2405.00651}.


\bibitem{milnor}
Milnor J.W., Stasheff J.D., Characteristic classes, \textit{Ann. of Math.
 Stud.}, Vol.~76, \href{https://doi.org/10.1515/9781400881826}{Princeton
 University Press}, Princeton, NJ, 1974.

\bibitem{oneill}
O'Neill B., The fundamental equations of a~submersion,
 \href{https://doi.org/10.1307/mmj/1028999604}{\textit{Michigan Math.~J.}}
 \textbf{13} (1966), 459--469.

\bibitem{palais}
Palais R.S., Homotopy theory of infinite dimensional manifolds,
 \href{https://doi.org/10.1016/0040-9383(66)90002-4}{\textit{Topology}}
 \textbf{5} (1966), 1--16.

\bibitem{Thurston}
Thurston W.P., Some simple examples of symplectic manifolds,
 \href{https://doi.org/10.2307/2041749}{\textit{Proc. Amer. Math. Soc.}}
 \textbf{55} (1976), 467--468.

\end{thebibliography}
\end{document}